\documentclass[11pt]{article}
\usepackage{amsmath}
\usepackage{amsfonts}
\usepackage{amssymb}
\usepackage{graphicx}
\newtheorem{thm}{Theorem}[section]
\newtheorem{cor}[thm]{Corollary}
\newtheorem{lem}[thm]{Lemma}
\newtheorem{prop}[thm]{Proposition}

\numberwithin{equation}{section}

\newcommand{\dx}{\,{\rm d}x}

\newcommand{\dt}{\,{\rm d}t}
\newcommand{\rd}{{\rm d}}

\def\LL{\mathrm{L}} 
\def\supp{\mathrm{supp}} 
\newcommand{\RR}{\mathbb{R}}

\newcommand{\bM}{\overline{M}}

\newcommand{\ve}{\varepsilon}
\newcommand{\vol}{{\rm Vol}}
\newcommand{\ren}{\RR^d}
\def\ee{\mathrm{e}} 
\def\dist{\mathrm{dist}} 
\def\qed{\,\unskip\kern 6pt \penalty 500
\raise -2pt\hbox{\vrule \vbox to8pt{\hrule width 6pt
\vfill\hrule}\vrule}\par}

\begin{document}
\title{\textbf{Positivity, local smoothing and Harnack
inequalities\\ for very fast diffusion equations}\\[7mm]
\normalsize{\textit{Dedicated to Luis Caffarelli for his upcoming
60\,$^{th}$ birthday}}}
\author{\Large Matteo Bonforte$^{\,a,\,b}$ 
~and~ Juan Luis V\'azquez$^{\,a,\,c}$\\} 
\date{} 

\maketitle

\

\begin{abstract}

We investigate qualitative properties of local solutions
$u(t,x)\ge 0$ to the fast diffusion equation, $\partial_t u
=\Delta (u^m)/m$ \ with $m<1$, corresponding to general
nonnegative initial data. Our main results are quantitative
positivity and boundedness estimates for locally defined
solutions in domains of the form  $[0,T]\times\RR^d$.
They combine into forms of new Harnack inequalities
that are typical of fast diffusion equations.
Such results are new for low $m$ in the so-called very fast
diffusion range, precisely for all $m\le m_c=(d-2)/d.$ The
boundedness statements are true even for $m\le 0$, while the
positivity ones cannot be true in that range.
\end{abstract}
\vspace{3cm}

\vspace{1cm}

\begin{itemize}
\item[(a)] Departamento de Matem\'{a}ticas, Universidad
Aut\'{o}noma de Madrid,\\ Campus de Cantoblanco, 28049 Madrid, Spain
\item[(b)] e-mail address:~matteo.bonforte@uam.es
\item[(c)] e-mail address:~juanluis.vazquez@uam.es
\end{itemize}


\newpage

\section*{\large Introduction}

We study qualitative and quantitative properties of  solutions
$u=u(t,x)$ of the nonlinear diffusion equation
\begin{equation}\label{FDE}
\partial_t u =\nabla\cdot(u^{m-1}\nabla u)=\Delta(u^m/m)
\end{equation}
in the whole parameter range $-\infty<m<1$, where it is called
Fast Diffusion Equation (FDE). We consider local nonnegative weak
solutions, defined in an open cylinder $Q$ of space-time
$\RR\times \RR^d$ with $d\ge 1$. Note that the factor $1/m$ in the
last expression is inessential when $m>0$ (up to a time rescaling,
$t'=t/m$) but becomes essential for $m< 0$, in order to obtain a
parabolic equation; for $m=0$ the last expression has to be
written as $\partial_t u =\Delta \log (u)$\footnote{We will always
interpret $u^m/m$ as $\log(u)$ when $m=0$. In the whole paper,
$\nabla$ indicates the gradient operator, $\nabla\cdot$ the
divergence operator, and $\Delta$ the Laplacian operator, all of
them taken with respect to the space variables, $x\in \RR^d$.}.

  Assuming the basic existence and uniqueness theory,
 \cite{DasKe},  \cite{VazLN}, we are interested in the qualitative
properties of the solutions such as boundedness, positivity, and
Harnack inequalities. For the FDE these properties depart from the
properties of the linear Heat Equation (case $m=1$), \cite{Widder},
and even more from the Porous Medium Equation (case $m>1$)\,,
\cite{VazBook}. Moreover, they are still partially understood when
$m$ is far from 1, precisely for $m\le m_c$ where $m_c=(d-2)/d$ \ is
called the {\sl first critical fast diffusion exponent}. Our goal
here is to obtain  bounds from above and below for the solutions in
that low range of exponents. We look for  precise quantitative
versions based on local estimates. Such estimates should be of
interest in developing a general theory of this equation in the
detail that is already known both for $m\ge 1$ and for $m_c<m<1$.

\subsection*{\large Precedents and problems}

The existence and uniqueness of weak solutions of the initial value
problem and other standard initial and boundary value problems for
the FDE, as well as  the main qualitative properties of the
solutions (such as the ones already mentioned, or the asymptotic
behaviour), are by now well understood when $m$ is close to one,
more precisely in the so-called {\sl good parameter range:} \
$m_c<m<1$.\footnote{With the extra restriction $m>0$ if $d=1$, the
case $-1<m<0$ and $d=1$ being somewhat different, cf. \cite{VazLN}.}
To be specific, when the problem is posed in the whole space, weak
solutions are uniquely determined by their initial data if $u_0$ is
a locally integrable nonnegative function, or even a locally finite
Radon measure. In that case, the solution is $C^\infty$ smooth and
positive for all $x\in \RR^d$ and $t>0$, and the initial data are
taken in the sense of initial trace, \cite{HP}, \cite{PiFDE},
\cite{DasKe}.  Solutions are bounded for data $u_0\in L^p(\RR^d)$
for any $p\ge 1$, and even for data in the Marcinkiewicz spaces
$M^p(\RR^d)$, $p> 1$, \cite{VazLN}. They are locally bounded under
the very mild restriction that $u_0$ is Radon measure, even if it is
not globally finite.

The theory of the FDE has been much less studied until recently in
the {\sl subcritical fast-diffusion range} $m<m_c$, even under the
condition $m>0$,  since essential difficulties have been found in
the different chapters of the theory, like existence, uniqueness,
and regularity. Note that $0<m<m_c$ is possible only if $d>2$. We
refer for background to the book  \cite{VazLN} that discusses in
some detail the range \ $m\le m_c$, even for $m\le 0$, along with
the cases $m>m_c$. Let us give an idea of the difficulties that
arise and that we address in our work below:

\medskip

\noindent  {\sc Boundedness.} Though weak solutions with data in
the spaces $\LL^p(\RR^d)$, $1\le p\le \infty$, exist and are unique
for $0<m<1$, counterexamples show that for $m<m_c$ these weak
solutions need not be bounded, and as a consequence they are not
smooth. The simplest such example seems to be the
separate-variables function
\begin{equation}\label{VSS}
U(t,x;T,x_0)=c\,\frac{(T-t)^{1/(1-m)}}{|x-x_0|^{2/(1-m)}}
\end{equation}
For every $m<m_c$, even $m\le 0$, there exists a suitable constant
$c(m,d)>0$ such that $U$ is a weak solution of the FDE in the
cylinder $Q=(0,T)\times \RR^d$, cf. \cite{VazLN}, page 80, but
obviously the solution never improves its initial regularity until
it extinguishes in finite time. The precise space regularity is
$U(\cdot,t)\in \LL^p_{loc}(\RR^d)$ for all $p<p_c$, where the
critical integrability exponent is  $p_c=d(1-m)/2$, which is larger
than 1 precisely for $m<m_c$, i.e., in the  subcritical range.

There is a positive result concerning boundedness, that is also tied
to the exponent $p_c$ : solutions with initial data in
$\LL^p(\RR^d)$ with $p>p_c$ become bounded and $C^\infty$ smooth for
all positive times as long as the solution does not disappear. This
{\sl smoothing effect} happens for all $p\ge 1$ if $m>m_c$, for
$p>1$ if $m=m_c$ (in the last  cases there is no problem of
disappearance). The results are sharp, cf. \cite{VazLN}.

\medskip

\noindent {\sc Extinction in finite time, EFT. } The above example
exhibits another typical feature of the Cauchy problem for
$m<m_c$, namely, the possible lack of positivity due to  EFT. The
occurrence of EFT depends on the type of problem we consider.

In the case of the Cauchy problem posed in $\RR^d$ with $d\ge 3$,
B\'enilan and Crandall  gave in \cite{BCr-cont} a proof of the
extinction in finite time, EFT,  of solutions of
 the FDE in the range $0<m<m_c$ when $u_0\in
\LL^p(\RR^d)$ with $p=p_c$. It is proved in \cite{VazLN} that EFT
occurs for the solutions with $m<m_c$ for all functions with initial
data in the Marcinkiewicz space $M^{p_c}(\RR^d)$, hence in
$\LL^{p_c}(\RR^d)$. We recall the EFT does not happen for the Cauchy
Problem when $m>m_c$.

 In the case of  the Cauchy-Dirichlet problem
posed in a bounded domain with zero boundary data, EFT happens for
all $0<m<1$. There is an interesting functional connection: we can
show that  EFT occurs if we have a global Poincar\'e and a Sobolev
inequality, and this result can be extended to more general
settings, such as Riemannian manifolds, as it has been done by the
authors in \cite{BGV-JEE}.  On the other hand, B\'enilan and
Crandall's proof  for the Cauchy problem is based only on the
Sobolev inequality, but it holds only in the lower range $m<m_c$.

\medskip

\noindent  {\sc Harnack inequalities.} Concerning finer regularity
properties, the possible occurrence of EFT is compatible with the
fact that nonnegative bounded solutions are positive, and
consequently $C^\infty$ smooth, as long as they are not identically
zero, i.e., before extinction. However, the existence of EFT for low
$m$ is tied to the breakdown of the standard forms of Harnack
inequalities, which are a strong tool in developing a regularity
theory. Obtaining some kind of Harnack inequality is therefore a
main research issue for $m\le m_c$ and has been an open problem for
some years. More specifically, we concentrate on parabolic lower
Harnack inequalities of the type called Aronson-Caffarelli estimates
 \cite{ArCaff}, and examine their consequences to obtain quantitative forms
of positivity. An extension work has been done in \cite{BV} for
$m_c<m<1$ but the method collapses for $m\le m_c$ due to the very
different properties of the solutions.
 As a consequence of our local
smoothing effect and of positivity estimates, we will obtain some
intrinsic Harnack inequalities of forward, elliptic or backward type,
which are new in this range.

In a recent preprint \cite{DGV}, DiBenedetto, Gianazza and Vespri
study the validity of intrinsic Harnack inequalities in the good
range $m>m_c$ and show, with an explicit counterexample, that any
kind of Harnack Inequality, intrinsic, elliptic, backward and
forward can not hold if $m<m_c$, for a fixed size of the intrinsic
cylinder, that is, if we fix the size of the parabolic cylinder
``a priori'' in terms of the value of $u$ at the center of the
cylinder $(t_0,x_0)$\,. They leave as an open problem to find
which kind of Harnack Inequalities, if any, are typical of the
very fast diffusion range $0<m<m_c$. In this paper we give an answer
to this intriguing problem.

\medskip

\noindent  {\sc Very singular range.} Most  the literature has
avoided the cases $m\le 0$, where the diffusivity $D(u)=u^{1-m}$
is very singular at $u=0$. Recently, it has been shown that a
large part of the theory of the subcritical range goes over to
this very singular range, on the condition of working with
solutions that ``are not too small''. See \cite{VazNonex} among
the older references, then \cite{DasDP}, and the books
\cite{DasKe}, \cite{VazLN} for a more complete reference. Note
that this recovers a subcritical range for dimensions $d=1,2$, and
also that we can study the interesting log-diffusion problems
where $m=0$, cf. \cite{Ham88},  \cite{Vld07} and the references.

More specifically, there is an extension of the results called
smoothing effects, whereby data in $\LL^p(\RR^d)$ with $p>p_c$ imply
bounded solutions for all $t>0$, and also the extinction in finite
time for data in $M^p(\RR^d)$, $p=p_c$. But a very different
situation happens for data in $\LL^p(\RR^d)$ with $1\le p<p_c$,
which is called {\sl immediate extinction}, whereby the solutions
obtained as limit of any reasonable approximation are identically
zero for all $t>0$. This makes it difficult to think of a general
study  of positivity. Immediate extinction happens for the
Cauchy-Dirichlet problem posed in a bounded domain with zero
boundary data for all $m\le 0$, $d\ge 1$. Our study of this range
is confined therefore to upper estimates.

\medskip

\noindent {\sc Comparison with elliptic  problems.} Part of the
difficulties of the FDE in the lower range of $m$ can be explained
by the intimate relation of the equation with  semilinear elliptic
theory.  This remarkable connection will be briefly explained in
Subsection \ref{Ell.Conn}.

\subsection*{ \large Results and organization}

Our work focuses on a peculiar feature of the FDE, which is the
existence of very strong local estimates. This was presumably first
mentioned in the paper by Herrero and Pierre \cite{HP}, 1985, who
get solutions in the whole range $0<m<1$ under the sole condition on
the initial data $u_0\in \LL^1_{loc}(\RR^d)$. Much of the subsequent
work has been influenced by the local character of the equation.
Here, we want push this idea to its final consequence concerning two
different areas: the question boundedness of local solutions, and
the question of positivity of nonnegative solutions, measured
quantitatively by so-called lower Harnack inequalities. We will then
combine the local upper and lower estimates, into a
full form of Harnack inequality.  While the boundedness
results hold for all $m<1$, positivity estimates are confined to
$0<m<1$ because of the possible occurrence of immediate extinction.
As we have said, the main interest of our results lies in their
application in the subcritical range, $m<m_c$. They are also new for
the critical exponent $m=m_c$.

Let us be more specific about the contents of the paper. It is
divided into  three main parts.

(I) The study of positivity and {\sl lower Harnack inequalities,}
both of local and global type. The first main contribution of the
paper is a {\sl parabolic lower Harnack inequality of the
Aronson-Caffarelli type} that is presented in Section \ref{sec-llb},
along with a detailed comparison  with the forms available for other
ranges of $m$. We devote Subsection \ref{sect.mdp} to prove the
lower estimate, Theorem \ref{Thm.Posit}, for a minimal problem. This
is extended in Subsection \ref{sect.pls} to  general solutions.
We then show that in the range  $m_c<m<1$
we can further eliminate the presence of the extinction time and recover
stronger estimates that are known in that range.
Subsection \ref{sect-loindT} discusses upper bounds for the
extinction time $T$ in terms of $\LL^p$ norms of the data, which give an alternative
type of lower  bound in the range where estimates depending only on $\LL^1$-norms of
the data are not true.

(II) The study of {\sl local upper bounds.} This takes two forms:
the first is the control of the evolution in  time of some spatial
$\LL^p_{\rm loc}$ norms, which is performed in Section
\ref{sect.lp}. Then, we get a local in space-time version of the
smoothing effect from $\LL^p_{\rm loc}$ into $\LL^\infty_{\rm
loc}$, an important regularity result that opens the road to
higher regularity and was known for $m>m_c$, and is false in
general for $m\le m_c$. We show in this paper that the estimate
holds $m<m_c$, on the condition that $p$ must be large enough.
We finally obtain the finest local upper estimates, called local
smoothing effects, in the form given in Theorem~\ref{thm.upper}\,,
just by combining the space-time smoothing effect and the
$\LL^p_{\rm loc}$ obtained in the first Section~\ref{sect.lp}.

 (III) Parabolic Harnack Inequalities. In Section \ref{sect.parabolicHI}, we combine
the local upper and lower estimates obtained in Parts I and II in
the form of parabolic Harnack inequalities of forward, backward and elliptic type,
together with an alternative form.

To conclude, we sketch a panorama of the obtained local estimates
depending on the ranges of $m$, together with general remarks, some
related open problem and a short review on related works. A final
Appendix contains some useful technical results.

\medskip

\noindent {\sc Notations.} We will work with weak solutions $u\ge
0$ of the FDE with $m<1$, defined in a cylinder $Q=\Omega\times
(T_0,T_1)$ for some domain $\Omega\subset \RR^d$ and $T_0<T_1$.
Usually, we take $T_0=0, T_1=T$. $T_1$ can be infinite and
$\Omega$ can be the whole space. In view of existing theory we may
assume that the solutions are positive and smooth as long as they
do not extinguish identically. We will be mostly interested in the
local theory where the space domain is bounded and the boundary
conditions are not taken into account. In deriving local estimates
it will be often sufficient to take as space domain a ball, which
we will denote by $B=B_R(x_0)$ or $B=B_{\lambda R}(x_0)$ for some
$\lambda>1$. We will frequently consider the annulus region
$A_{R,\lambda}=B_{\lambda R}\setminus B_R$.  As indicated before,
we put
$$
m_c=\frac{d-2}d, \qquad p_c=\frac{d(1-m)}{2}\,.
$$
We have pointed out that $p_c>1$ if and only if $m<m_c$. We will
take integrability exponents $p\ge 1$ if $m>m_c$, $p>p_c$ if $m\le
m_c$. Moreover, for $p\ne p_c$ we set
\begin{equation}
\vartheta_p=\frac{1}{2p-d(1-m)}\,,
\end{equation}
which is positive if and only if $p>p_c$.

\section{\large  \textsc{Part I.} Local lower bounds}
\label{sec-llb}

The first part of the paper addresses the question of quantitative
estimates of positivity. The exponent range in this part is $0<m<1$,
since it is well known that the FDE does not admit solutions of the
Dirichlet  problem with zero boundary data when $m\le 0$, thus
blocking any possibility of a general local positivity theory in
that range \cite{VazNonex, VazLN}. Our main contribution is a
parabolic inequality in the spirit of the one obtained by Aronson
and Caffarelli \cite{ArCaff} in their path-breaking paper for $m>1$,
and the ones produced by the authors in \cite{BV} for $m_c<m<1$. The
purpose of such formulas is giving quantitative information on the
positivity of solutions at later times in terms of information on
$\LL^p$ norms of $u$ at a former time that we take as $t=0$. This is
why they are called parabolic lower Harnack formulas.

We take  $0<m<1$ and consider a  $u$ be a  local, nonnegative weak
solution of the FDE defined in a cylinder $Q=(0,T)\times \Omega$,
taking initial data $ u(0,x)=u_0(x)$ in $ \Omega$ and having
finite extinction time $T$. We make no assumption on the boundary
condition (apart from nonnegativity). For ease of proof we will
assume that the solutions are smooth so that the different
computations and comparison results are valid. This assumption is
then eliminated by approximation, which is justified according to
known theory.

\begin{thm}\label{Thm.Posit}
Let $0<m<1$ and let $u$ be the solution to the FDE under the above
assumptions. Let $x_0$ be a point in $\Omega$ and let
$d(x_0,\partial \Omega)\ge 3R$. Then the following inequality
holds for all $0<t<T$
\begin{equation}\label{AC.Estimates}
R^{-d}\int_{B_R(x_0)}u_0(x)\dx
  \le
  C_1\,R^{-2/(1-m)}\,t^{\frac{1}{1-m}}+
  C_2\,T^{\frac{1}{1-m}}R^{-2}\;t^{-\frac{m}{1-m}}\;u^m(t,x_0).
\end{equation}
with $C_1$ and $C_2$ given  positive constants depending only on
$d$. This implies that there exists a time $t_*$ such that
 for all $t\in(0,t_*]$
\begin{equation}\label{Flux.Posit.Sub.Dirichlet}
 u^m(t,x_0)\ge C_1'\,R^{2-d}\|u_0\|_{\LL^1({B_R})} T^{-\frac{1}{1-m}}\;
t^{\frac{m}{1-m}}.
\end{equation}
where $C_1'>0$ depends only on $d$; $t_*$ depends on $ R$ and
$\|u_0(x)\|_{\LL^1({B_R})}$ but not on $T$. \end{thm}

\noindent {\bf Simplified version.} The dependence on the parameters
makes the formula apparently complicated. But it can be reduced to a
simpler, equivalent one. Actually, we may assume that $x_0=0$ by
translation. Given $R>0$ and  $M=\int_{B_R(0)}u_0(x)\,dx>0$, we use
the rescaling
\begin{equation}
{u}(t,x)=\frac{M}{R^d}\,\widehat{u}\left(\frac{t}{\tau}, \frac{x}{R}
\right)\,,\qquad \tau= R^{2-d(1-m)}M^{1-m},
\end{equation}
to pass from a solution with mass $M$ in the ball of radius
$R$ to a solution $\widehat{u}$ with mass 1 in the ball of radius 1.
So we only need to prove the version with $M=R=1$ to get the full
version. The scaling is simpler for $m=m_c$ where $\tau=M^{1-m}$. Of
course, the extinction time has to be rescaled accordingly,
$T=R^{2-d(1-m)}M^{1-m}\widehat{T}$.

\medskip

\noindent {\bf Improvements.} As stated, estimate
\eqref{AC.Estimates} applies only to solutions with finite
extinction time, and it involves  the value of the extinction time
$T$ in an explicit way; both things can make it impractical.
However, a simple comparison argument shows that we only need to
estimate from below any subsolution. In particular, we may replace
the solution under consideration by the solution of the problem with
initial data $u_0(x)\chi_{B_R(x_0)}(x)$, and zero Dirichlet boundary
conditions on $x\in
\partial B_{3R}(x_0)$. Let us call this problem  \sl minimal problem\rm \
for the given data. The extinction time of the corresponding
solution will be called the {\sl minimal life time} \ of such domain
and data,  $T_m(u_0,B)$. Clearly, $T_m(u_0,B)\le T (u)$.

\begin{cor} \label{coroll1.2} \sl
The above positivity result holds with $T(u)$  replaced by the
minimal life time \ $T_m(u_0,B)$, $u$ is defined in $Q_{T}$,  and
the estimate applies for $0<t<T'$ with $T'=\min\{T, T_m\}$.
\end{cor}

This modified result is specially interesting in the range $1>m>m_c$
where the solutions of the Cauchy Problem do not vanish. On the
other hand, it is known that $T_m$ is finite if $u_0$ satisfies some
local integrability conditions \cite{DiDi, VazBook}.
\subsubsection*{\large Comparison with the estimate for the PME and other FDE}

\noindent {\sc The PME.}  Let us write Aronson-Caffarelli's result
\cite{ArCaff} for $m>1$ with a similar notation:
\begin{equation}\label{ACineq}
R^{-d}\int_{B_R(x_0)}u_0(x)\dx   \le
  C_1\,R^{2/(m-1)}\,t^{-\frac{1}{m-1}}+C_2\,R^{-d}t^{d/2}u^{1+ (d(m-1)/2)}(t,x_{0}).
\end{equation}
We recall that this formula is valid for all nonnegative weak
solutions of the PME defined in the whole space.  The form of the
first term in the right-hand side is the same in both results,
\eqref{AC.Estimates} and \eqref{ACineq}. This term plays the role of
blocking the positivity information when it is large relative to the
left-hand side integral, and allowing for such information when it
is small. The critical time at which we begin to get positivity
information is obtained by making this term a fraction of the
left-hand side, i.e., for
\begin{equation}\label{crit.time}
t_c=c(m,d)\|u_0\|^{1-m}_{\LL^1(B_R(x_0))}R^{2+d(m-1)}.
\end{equation}
But since the exponents have just the opposite sign in the above
expressions  for $m>1$ and $m<1$, the consequences are qualitatively
very different: the information on positivity happens for us when
$t$ is smaller than $t_*$, while for the PME it happens when $t$ is
larger. This is in accord with the basic properties of these
equations, which the present inequalities faithfully reproduce.
Rescaling allows to check the inequality only at $t=1$ for $R=1$,
and in that case we only have to prove that there are constants
$M_{0}=M_{0}(n,m)$ and $k=k(m,d)$ such that for $M\geqslant M_{0}$
\begin{equation}\label{Harn.estim.tildeu}
u(0,1)\geqslant k\,M^{2/(d(m-1)+2)}.
\end{equation}
As to the  second term, it is different. We cannot expect to have
the A-C term in the range $m<m_c$ since then the exponent of $u$
would be negative. In fact, the proof of  \cite{ArCaff} uses
conservation of mass that is not valid for the fast diffusion
equation in the low $m$ range.

\medskip

\noindent{\sc The good FDE.} The validity of the Aronson-Caffarelli
formula was extended by the authors in \cite{BV} to local solutions
of the FDE in the good exponent range $m_c<m<1$, and the already
mentioned sign change in the exponents implies that we get good
lower estimates for $0<t\le t_*$. Moreover, we can continue these estimates thanks to
the fortunate circumstance that we have further differential
inequalities, like $\partial_t u\ge -Cu/t$ in the case of the Cauchy
problem, which allow for a continuation of the lower bounds for
$t\ge t_*$ with optimal decay rates in time. The final form is
\begin{equation}\label{Posit.FDE}
 u(t,x)\ge \bM_{R}(x_0)\, H(t/t_c), \quad
 \bM_{R}(x_0)=R^{d}\int_{B_R(x_0)}u_0\,dx.
\end{equation}
The critical time is defined as in \eqref{crit.time}; the function
$H(\eta)$ is defined as $K \eta^{1/(1-m)}$ for $\eta\le 1$ while
$H(\eta)=K \eta^{-d\vartheta} $  for $\eta\ge 1$, with $K=K(m,d)$.
Note that for $0<t<t_c$ the lower bound means $u(t,x_0)\ge
k(m,d)(t/R^2)^{1/(1-m)}$ which is independent of the initial mass.

 \medskip

 \noindent{\sc Eliminating the  time $T$.}
A natural question is to try to recover this sharp results of the
good fast diffusion range via the present methods. If one wants to
do that, one needs upper estimates for the minimal life time, that
is upper estimates for the extinction time for the MDP, in terms of
the $\LL^1$-norm on the ball $B_{R_0}$.  We prove the following
result.

\begin{thm}\label{GoodFDE.Posit} Let $m_c<m<1. $ Then,
{\rm (i)} We have sharp upper and lower estimates for the extinction
time for the Dirichlet problem on any ball $B_R$ of the form:
\begin{equation}\label{GoodFDE.sharp.ext.est}
c_1\|u_0\|_{\LL^1(B_{R/3})}^{1-m}R^{2-d(1-m)}\le T \le c_2
\|u_0\|_{\LL^1(B_R)}^{1-m}R^{2-d(1-m)}.
\end{equation}

\noindent {\rm (ii)} In that range of $m$ the lower estimates of
Theorem \ref{Thm.Posit} imply the lower Harnack inequalities of {\rm
\cite{BV, DKV, DiBenedettourbVesp, DGV},} in the form
\begin{equation}\label{GoodFDE.Low.Harnack}
u(t,x_0)\ge c_{m,d}\left[\frac{t}{R^{2}}\right]^{\frac{1}{1-m}}
\end{equation}
for any $0<t<t^*$ and any $x\in B_R$, where $t_*$ is given by
\eqref{form.t_*}.
\end{thm}

This result shows that the form of the lower bounds given in Theorem
\ref{Thm.Posit} is sharp, since it allows to obtain sharp local
lower bounds not only in the good fast diffusion range. And it also
applies in the very fast diffusion range, that is the new
interesting part of this paper. We are thus led to the question  of
eliminating all extinction times from the estimate, i.e., replacing
$T$ or $T_m$ by some information on the initial data, also in the
range $0<m<m_c$.

\begin{thm}\label{Thm.Posit.v2}
Let $0<m<m_c$ and let $u$ be the solution to the FDE under the above
assumption that $u_0\in  \LL^{p_c}_{loc}(\RR^d)$. Let $x_0$ be a point
in $\Omega$ and let $d(x_0,\partial \Omega)\ge 3R$. Then, the
following inequality holds for all $0<t<T$
\begin{equation}\label{AC.Estimates.v2}
R^{-d}\|u_0\|_{\LL^{1}(B_R(x_0))}
  \le
  C_1\,R^{-2/(1-m)}\,t^{\frac{1}{1-m}}+
  C_3\,\|u_0\|_{\LL^{p_c}(B_R(x_0))} R^{-2}\;t^{-\frac{m}{1-m}}\;u^m(t,x_0).
\end{equation}
with $C_1$ and $C_3$ given positive constants depending on $d$.
\end{thm}

 We can also obtain formulas in terms of the norms
$\|u_0\|_{\LL^p(B_R(x_0)}$ for all $p>p_c$, that can be seen below.

\subsubsection*{\large Obstruction to a simpler  estimate with $\LL^1$ norm}

The presence of the extinction time $T$ in the lower estimates, or
equivalently of some $\LL^p$ norm of the initial data, is a drawback
in the formulas that is not present in the original
Aronson-Caffarelli estimate for $m>1$, or in the version of the
authors for $m\in (m_c,1)$ in the whole space. But it is a
consequence of the `bad' behaviour of the fast diffusion equation
for low values of $m$, a fact that can be seen in different ways.

Thus, we will show here that the local lower estimates cannot depend
only on the local $\LL^1$ norm of the initial data  when $0<m\le
m_c$. We do it  by means of a counterexample based on the behaviour
of solutions with data that approximate a Dirac delta. We solve the
FDE for smooth and positive initial data $\varphi(x)\in
\LL^1(\RR^d)$ with integral equal to 1. We assume  that $\varphi$ is
radially symmetric, compactly supported and decreasing with $|x|$.
We obtain a smooth and positive solution $u(t,x)$ defined in a
cylinder $Q_{T_1}$ and vanishing identically at some $t=T_1$. The
scale invariance of the equation implies that the solution
corresponding to data $\varphi_k(x)=k^{d}\varphi(kx)$ is
\begin{equation}
u_k(x)=k^d u(k^{-\sigma}t,k x), \quad \sigma=d(1-m)-2>0,
\end{equation}
so that it has extinction time $T_k=T_1\,k^{\sigma}$. As
$k\to\infty$ it is clear that $u_k(0,t)$ converges to the Dirac
delta. We also observe that $T_k\to\infty$, so that we lose the
previous estimates. On the other hand, we see that losing the
estimates is inevitable. If we consider a point $x_0$ very close to
$x=0$ and take a radius $R>|x_0|$, then
$\|u_k(0,x)\|_{\LL^1(B_R(x_0))}=1$. However, by continuity of $u$
with respect to the initial data at $t=0$, $x$ large, we have
$$
u_k(t,x_0)=k^d u(k^{-\sigma}t,kx_0)\to 0
$$
(note that $\int_{\ren }u_k(t,x) \,dx\le 1$ at all times). This
means that no lower estimate could be uniformly valid for this
sequence.

A scaling argument was used by Brezis and Friedman \cite{BrFr} to
prove that there exist no weak solutions with initial data a Dirac
delta.

\subsection{\large Positivity for a ``minimal'' Dirichlet Problem}
\label{sect.mdp}

We will assume that $0<m<1$ in the study  of positivity (cf. the
comment in the Introduction). Since $m>0$  we eliminate the factor
$1/m$ from  equation \eqref{FDE} for simplicity without loss of
generality.  As a preliminary step, we first prove positivity for
a problem posed on a ball of radius $R_0$, zero boundary data and
particular initial data. Since the problem of getting quantitative
positivity estimates has been successfully studied in \cite{BV} in
the range $m_c<m<1$, the techniques we introduce are mainly aimed
at producing positivity in the cases $0<m\le m_c$, where previous
methods failed.

Specifically,  we shall consider the following Dirichlet problem on
the ball $B_{R_0}\subset \RR^d$:
\begin{equation}\label{Sub.Dirichlet.Problem}
\left\{\begin{array}{lll}
 \partial_t u =\Delta (u^m) & ~ {\rm in}~ Q_{T,R_0}=(0,T)\times B_{R_0}\\[2mm]
 u(0,x)=u_0(x) & ~{\rm in}~ B_{R_0},~~~{\rm and}~~~ \supp(u_0)\subseteq B_{R}\\[2mm]
 u(t,x)=0 & ~{\rm for}~  t >0 ~{\rm and}~ x\in\partial B_{R_0}\,,
\end{array}\right.
\end{equation}
where  $R_0>2R>0$. We only consider
nonnegative data and solutions. The problem admits a unique
solution $u\in C([0,\infty):  \LL^1(B_{R_0})) $ for every $
u_0\in\LL^{1} \big(B_{R_0}\big), $ \cite{BCr}. We will refer to
this problem as the minimal Dirichlet problem, or more briefly,
the {\sl minimal problem}, because obtaining positivity for
solutions to this problem implies in an easy way local positivity
for any other problem, thanks to the comparison principle.   The
solution vanishes in finite time; let $T>0$ be the finite
extinction time, shortly FET. Later on we would like to eliminate
the dependence of the results on $T$ and make the estimates depend
only on the initial data, see Section \ref{sect-loindT}.

Our goal is to obtain  positivity with a quantitative estimate for
this ``minimal'' problem. Our most novel idea consists in passing
the information on the initial data via the flux of the solution on
the boundary of the ball $B_{2R}$ into an averaged positivity result
outside the ball for times that are not too small, more precisely on
the annulus $A_0:=B_{R_0}\setminus B_{2R}$. This property can be
interpreted as the expansion of positivity outside a ball in which
the initial datum has nonzero mean. It is in some sense it is analogous
to the expansion of positivity already introduced by DiBenedetto et
al.\,, see e.g. \cite{DK, DGV, DiBenedettourbVesp} for the upper $m$-range\,. The expansion of
positivity turns out to be a key tool in proving lower Harnack also
in our case.

Once we have proved that positivity spreads out from a ball, then
for suitable positive times the mean value of the solution on an
annulus is positive. We then "fill the hole" in the middle using
 Aleksandrov's Reflection Principle, cf. the Appendix and
\cite{BV}. In this way we arrive at the positivity result in the
inner ball for any positive time.

\subsubsection{Flux and transfer of positivity}

We start the proof of the positivity results for the minimal problem
by a result on mass transfer to an outside annulus  based on the
flux across an internal boundary. We recall that $R_0>2R$ and
$A_0:=B_{R_0}(x_0)\setminus B_{2R }(x_0)$. In order to simplify the
final formulas, we write $\lambda=R_0/2R>1$ (we take for instance
$R_0=3R$).

\begin{lem}[Flux Lemma]\label{Flux.Lemma}
If $u$ is a positive smooth solution of the Minimal Problem {\rm
\eqref{Sub.Dirichlet.Problem}} in $Q_T$ with extinction time $T>0$.
Then, the following estimate holds true
\begin{equation}\label{Flux.Lemma.Posit}
k_0\,(R_0-2R)^2 \int_{B_{R_0}}u(s,x)\,dx \le
\int_s^T\int_{A_{0}}u^m\dx \dt,
\end{equation}
for any $0\le s \le T$, and any $0<2R<R_0$, and for a suitable
constant $k_0=k_0(d)$.
\end{lem}

\noindent {\sc Proof.~}  We shall use a $C^\infty$ test function
$\varphi(x)$ that is supported in the ball $B_{R_0}$ and takes the
value 1 in $B_{2R}.$  It is clear that we can choose $\varphi$
such that there exist a constant $k_0 >0$ depending only on $d$
such that
\begin{equation}\label{GTF.Est}
\big|\Delta \varphi(x)\big|
 \le \frac{k_0^{-1}}{(R_0-2R)^2}\;,
\end{equation}
 Let $0\le s < t\le T$. We compute
\begin{equation*}
\begin{split}
\int_s^t\int_{A_{0}}\partial_t u \, \varphi\dx \dt
 &=\int_s^t\int_{A_{0}}\Delta(u^m) \varphi \dx \dt
 =\int_s^t\int_{A_{0}} u^m \Delta\varphi \dx\dt
  +\int_s^t\int_{\partial B_{2R}}\partial_\nu (u^m)\varphi \rd\sigma \dt\\
 &+\int_s^t\int_{\partial B_{R_0}}\big[\partial_\nu (u^m)\varphi -
 u^m \partial_\nu \varphi \big]\rd\sigma \dt
  -\int_s^t\int_{\partial B_{2R}} u^m \partial_\nu \varphi
  \rd\sigma \dt.
\end{split}
\end{equation*}
We remark that the last three integrals vanish since $\varphi$ and
$u\equiv $  vanishes identically near the boundary  $\partial
B_{R_0}$, and $\partial_\nu\varphi\equiv 0$ on $\partial B_{2R}$.
We also  have
\[
\int_s^t\int_{A_{0}}\partial_t u \, \varphi \dx \dt
 =\int_{A_{0}}u(t,\cdot) \varphi \dx
 -\int_{A_{0}}u(s,\cdot) \varphi \dx
\]
\noindent Hence,
\begin{equation*}
\int_{A_{0}}u(t,\cdot) \varphi \dx
 -\int_{A_{0}}u(s,\cdot) \varphi \dx
 =\int_s^t\int_{A_{0}} u^m \Delta\varphi \dx \dt
 +\int_s^t\int_{\partial B_{2R}}\partial_\nu (u^m)\varphi \rd\sigma\dt.
\end{equation*}
We will use this equality with $t=T$, $T=T(u_0)$ being the finite
extinction time for the solution to Problem
\eqref{Sub.Dirichlet.Problem}, so that we obtain
\begin{equation}\label{Flux.Eq.1}
 \int_{A_{0}}u(s,\cdot) \varphi \dx
 =-\int_s^T\int_{A_{0}} u^m \Delta\varphi \dx \dt
 +\int_s^T\int_{\partial B_{2R}}\partial_{\nu^*} (u^m)\varphi \rd\sigma\dt\,,
\end{equation}
where $\nu^*$ is the exterior normal to $B_{2R}$, which is the
opposite of $\nu$ which is the exterior normal to the inner
boundary of $A_{0}$, so that
$\partial_{\nu^*}(u^m)=-\partial_{\nu}(u^m)$.

On the other hand, a simple calculation shows that
\begin{equation*}
\int_{B_{2R}} (u(t,x)-u(s,x))\dx
 =\int_s^t\int_{B_{2R}}\partial_t u \dx \dt
 =\int_s^t\int_{B_{2R}}\Delta (u^m)\dx \dt
 =\int_s^t\int_{\partial B_{2R}}\partial_{\nu^*}(u^m)\rd\sigma\dt.
\end{equation*}
Letting $t=T$, with $T$ as above, we obtain
\begin{equation}\label{Flux.Eq.2}
-\int_{B_{2R}}u(s,x)\dx
 = \int_s^T\int_{\partial B_{2R}}\partial_{\nu^*}(u^m)\rd\sigma\dt.
\end{equation}
\noindent Joining equalities \eqref{Flux.Eq.1} and
\eqref{Flux.Eq.2} we get
\begin{equation*}
\int_{B_{R_0}}u(s,x)\dx=\int_{B_{2R}}u(s,x)\dx+\int_{A_0}u(s,x)\dx
 =-\int_s^T\int_{A_0}u^m\,\Delta\varphi\dx \dt
\end{equation*}
We conclude by using estimates \eqref{GTF.Est} for
$\Delta\varphi$: for any $0<2R<R_0$ we then get
\begin{equation}
\int_{B_{R_1}}u(s,x)\dx\le\int_{B_{R_0}}u(s,x)\dx=-\int_s^T\int_{A_0}u^m\,\Delta\varphi\,\dx \dt
     \le \dfrac{k_{0}^{-1}}{(R_0-2R)^2} \int_s^T\int_{A_{0}}u^m\,\dx \dt.
\end{equation}
The proof is  complete.\qed

\medskip

\noindent\textbf{Remark. }\textsl{Lower Bound on the Extinction
Time.} As a first consequence of this Lemma we can easily obtain
useful lower estimates for the FET:
\begin{equation*}
\begin{split}
k_0\,(R_0-2R)^2 \int_{B_{R_0}}u(s,x)\,&\dx
  \le \int_s^T\int_{A_{0}}u^m\dx \dt
  \le (T-s)\;\vol(A_0)\;\int_{A_{0}}u^m(\overline{s},x)\frac{\dx}{\vol(A_0)}\\
 &\le (T-s)\;\vol(A_0)\;\left[\int_{B_{R_0}}u(\overline{s},x)\frac{\dx}{\vol(A_0)}\right]^m\\
 &\le (T-s)\;\vol(A_0)^{1-m}\;\left[\int_{B_{R_0}}u(s,x)\dx\right]^m\\
\end{split}
\end{equation*}
where in the first step we have used the mean value theorem for the
time integral (see details in Step 2 of next section), with
$\overline{s}\in(s,T)$, in the second step the H\"older inequality,
and in the third step we used the contractivity of the global
$\LL^1(B_{R_0})$-norm. Letting then $s=0$ gives the desired lower
bound, once we notice that $\supp(u_0)\subseteq B_R$
\begin{equation}
k_0\,(R_0-2R)^2\left[\frac{\int_{B_{R}}u_0 \dx}{\vol\big(A_{0}\big)}
\right]^{1-m}\le T\,.
\end{equation}

\subsubsection{Pointwise lower estimate for initial
times}\label{Sec.Pos._Init}

 We have just shown how positivity of the initial datum
propagates on the annulus in the weak form of a positive space-time
mean value. We will now see that this is sufficient to fill the hole
inside the annulus. As in the study of the exponent range $m_c< m<1$
performed in \cite{BV}, the estimate uses a {\sl critical time} that
is defined in terms of the initial norms. In the present case it is
given by
\begin{equation}\label{t*.prima.forma}
 t_*\,:=
  \frac{k_0}2\,(R_0-2R)^2\left[\frac{\int_{B_{R}}u_0 \dx}{\vol\big(A_{0}\big)} \right]^{1-m}
\end{equation}
where $k_0$  as in the Flux Lemma \ref{Flux.Lemma}. Note in passing
that the positivity result that follows, formula
\eqref{Flux.Posit.Center}, implies that this quantity is less than
$T$.

\noindent Obtaining the lower bound needs several  steps.

\noindent$\bullet~$\textsc{Step 1.} \textsl{ Time Integrals.}
H\"older's inequality, together with the fact that the global
$\LL^1\big(B_{R_0}\big)$-norm decreases, gives
\[\begin{split}
\int_{A_{0}}u(t,x)^m\dx
 &\le\vol\big(A_{0}\big)^{1-m}\left[\int_{A_{0}}u(t,x)\dx\right]^m
  \le\vol\big(A_{0}\big)^{1-m}\left[\int_{B_{R_0}}u(t,x)\dx\right]^m\\
 &\le\vol\big(A_{0}\big)^{1-m}\left[\int_{B_{R_0}}u(0,x)\dx\right]^m
  =\vol\big(A_{0}\big)^{1-m}\left[\int_{B_{R}}u_0\dx\right]^m
\end{split}
\]
since $\supp(u_0)\subseteq B_R$\,. For any $0\le s\le t$ we then
have
\[
\int_s^t\int_{A_{0}}u(\tau,x)^m\dx \rd\tau
 \le\vol\big(A_{0}\big)^{1-m}\;\left[\int_{B_{R}}u_0\dx\right]^m\;(t-s)
\]
We use this estimate together with estimate \eqref{Flux.Lemma.Posit}
to get
\begin{equation}\label{Flux.Eq.3}\begin{split}
k_0\,(R_0-2R)^2\int_{B_R}u_0(x)\dx
 & \le \int_0^T\int_{A_{0}}u^m\dx \dt
   = \int_0^{t_*}\int_{A_0}u^m\dx \dt
   + \int_{t_*}^T\int_{A_{0}}u^m\dx \dt\\
 & \le
 \vol\big(A_{0}\big)^{1-m}\;\left[\int_{B_{R}}u_0\dx\right]^m\;t_*
   + \int_{t_*}^T\int_{A_{0}}u^m\dx \dt\\
\end{split}
\end{equation}
In view of the definition of $t_*$ we can eliminate one term and get
\begin{equation}
k_2\,(R_0-2R)^2\,\int_{B_R}u_0(x)\dx\,
 \le \int_{t_*}^T\int_{A_{0}}u^m\dx \dt\;.
\end{equation}
with $k_2=k_0/2$. In particular, this means that the left-hand side
remains strictly positive.

\medskip

\noindent$\bullet~$\textsc{Step 2.}  We introduce the function
\[
Y(t)=\int_{A_{0}}u^m(t,x)\dx,
\]
and apply the mean value theorem -for the time integral- to
prove that there exists $t_1\in[t_*,T ]$ such that $\int_{t_*}^{T}
Y(t)\dt
 = (T-t_*)\; Y\big({t_1}\big).$ In other words,
\begin{equation*}
\int_{t_*}^T\int_{A_{0}}u^m(t,x)\dx \dt
 =\big(T-t_*\big)\int_{A_{0}}u^m\big(t_1,x\big)\dx.
\end{equation*}
Using now the estimate obtained in the previous step, we conclude
that there exists $t_1\in \big[t_*\,,\,T)$, such that
\begin{equation*}
k_2\,(R_0-2R)^2\int_{B_R}u_0(x)\dx \le
\int_{t_*}^T\int_{A_{0}}u^m(t,x)\dx \dt
=\big(T-t_*\big)\int_{A_{0}}u^m\big(t_1,x\big)\dx,
\end{equation*}
and this  implies that for some $t_1\in [t_*,T]$ we have
\begin{equation}\label{Flux.Eq.4}
\frac{k_2\,(R_0-2R)^2}{T} \int_{B_R}u_0(x)\dx
 \le \frac{k_2\,(R_0-2R)^2}{\big(T-t_*\big)} \int_{B_R}u_0(x)\dx
 \le \int_{A_{0}}u^m\big(t_1,x\big)\dx.
\end{equation}

\noindent$\bullet~$\textsc{Step 3.} \textsl{ Aleksandrov Principle.
Positivity at the critical time.} We can now use the Aleksandrov
Principle to deduce positivity at $x_0$ from  inequality
\eqref{Flux.Eq.4}. In fact,
\begin{equation}\label{Flux.Aleks}
\int_{A_{0}}u^m\big(t_1,x\big)\dx\le \vol(A_{0})\,u^m(t_1,x_0)
\end{equation}
where $x_0\in\RR^d$ is the center of the ball $B_{R_0}$, since we
know from Aleksandrov Principle, that $u(t,x)\le u(t,x_0)$ for any
$x\in A_{0}$ and any $t>0$ (see Appendix for details).

\noindent Joining inequality \eqref{Flux.Eq.4} and
\eqref{Flux.Aleks} we obtained that  there exists a
$t_1\in\big[t_*\,,\,T)$ such that
\begin{equation}\label{Flux.Posit.tc}
\frac{k_2\,(R_0-2R)^2}{\vol(A_0)\,T}\;\int_{B_R}u_0(x)\dx
  \le \,u^m(t_1,x_0)
\end{equation}

\noindent$\bullet~$\textsc{Step 4.} \textsl{ Positivity backward in
time.} The last step consists in obtaining a lower estimate when
$0\leq t\le t_1$. This argument is  based on B\'enilan-Crandall's
differential estimate, cf. \cite{BCr}\; :
\[
\partial_t u (t,x)\leq\frac{u(t,x)}{(1-m)t}
\]
that is valid for all nonnegative solutions of this initial and
boundary value problem. It easily implies that the function:
\[
u(t,x)t^{-\frac{1}{1-m}}
\]
is non-increasing in time, thus for any $t\in (0,t_1]$ we have
that
\[
u(t_1,x)\;
 \le\; t^{-\frac{1}{1-m}}\;t_1^{\frac{1}{1-m}}u(t,x)\;
 \le\; t^{-\frac{1}{1-m}}\;T^{\frac{1}{1-m}}\;u(t,x)\;
\]
since $t_1\le T$\,. It is now sufficient to apply inequality
\eqref{Flux.Posit.tc} to the l.h.s. in the above inequality to get
\begin{equation}\label{Flux.Posit.Center}
\frac{k_2\,(R_0-2R)^2}{\vol(A_0)\,T}\;\int_{B_R}u_0(x)\dx
  \le \,u^m(t_1,x_0)
  \le\;t^{-\frac{m}{1-m}}\;T^{\frac{m}{1-m}}\;u^m(t,x_0)
\end{equation}
This is the  inequality we were looking for.

\noindent $\bullet~$\textsc{Step 5.} In order to simplify the
final formulas, it is convenient to make a choice for the ratio
$\lambda=R_0/2R>1$ (for instance  $R_0=3R$). The  formula for
$t_*$ becomes
\begin{equation}\label{form.t_*}
t_*= c_0'\,R^{2-d(1-m)}\, \|u_0 \|_{\LL^1(B_{R})}^{1-m}
\end{equation}
and $c_0'>0$ depends only on $d$ and $\lambda$. We have proved the
following positivity result.

\begin{thm}\label{Flux.Thm.Posit.2}
Let $0<m<1$ and let $u$ be the solution to the Minimal Problem
{\rm {\rm \eqref{Sub.Dirichlet.Problem}}} and let $T=T(u_0)$ be
the MET. Then $T\ge  2t_*$, and the following inequality holds
true for all $t\in(0,t_*]$
\begin{equation}\label{Flux.Posit.Sub.Dirichlet.2}
 u^m(t,x_0)\ge c_1'\,R^{2-d}\|u_0\|_{\LL^1({B_R})} T^{-\frac{1}{1-m}}\;
t^{\frac{m}{1-m}}.
\end{equation}
where $c_1'>0$ depends only on $d$.
\end{thm}

For the particular time $t=t_*$ we get
$$
 \left( \frac{u(t_*,x_0)}{\oint_{B_R} u(0,x)\,dx}\right)^{m}
 \ge c'_2\,(R^2/T)^{1/(1-m)}
\oint_{B_R} u(0,x)\,dx.
$$

\subsubsection{Estimate of Aronson-Caffarelli type for Very
Fast Diffusion}

It is interesting to present the above result in the form that has
been used by Aronson and Caffarelli in their work \cite{ArCaff}. We
have to argue as follows: we have arrived at the following
alternative
\[
\mbox{either}\qquad t^*<t\qquad\mbox{or}\qquad
 \frac{k_2(\lambda^2-1)}{\omega_d\,(\lambda^d-1)R^{d-2}\,T^{\frac{1}{1-m}}}\;\int_{B_R}u_0(x)\dx
  \le\;t^{-\frac{m}{1-m}}\;u^m(t,x_0)
\]
Writing the expression of $t^*$, we either  have
\[
R^{-d}\int_{B_R}u_0(x)\dx\le\,C_1\,R^{-2/(1-m)}\,t^{\frac{1}{1-m}},
\quad C_1=\frac{\omega_d(\lambda^d-1)}
{k_1^{\frac{1}{1-m}}\,(\lambda-2)^{\frac{2}{1-m}}},
\]
or
\[
R^{-d}\int_{B_R}u_0(x)\dx
  \le C_2\,T^{\frac{1}{1-m}}R^{-2}\;t^{-\frac{m}{1-m}}\;u^m(t,x_0),
  \quad C_2=\frac{\omega_d\,(\lambda^d-1)}{k_2\,(\lambda^2-1)}
\]
We now sum up the two expressions to get
\begin{equation}\label{AC.Estimates1}
R^{-d}\int_{B_R}u_0(x)\dx
  \le
  C_1\,R^{-2/(1-m)}\,t^{\frac{1}{1-m}}+C_2\,T^{\frac{1}{1-m}}R^{-2}\;t^{-\frac{m}{1-m}}\;u^m(t,x_0).
\end{equation}
with $C_1$ and $C_2$ given constants depending on $d$ and
$\lambda>2$.

\subsection{\large Proof of Theorem \ref{Thm.Posit} and Corollary \ref{coroll1.2} }\label{sect.pls}

The previous results will now be used to prove uniform positivity on
balls for any local solution of the FDE problem. We proceed as
follows: let $u$ be a positive and continuous weak solution of the
FDE defined in $Q=(0,T)\times \Omega$ taking initial data $
u(0,x)=u_0(x)$ in $ \Omega$. We make no assumption on the boundary
condition (apart from nonnegativity). Let us select a point
$x_0\in\Omega\subseteq\RR^d$. Select two radii $R_0\ge 3R>0$ so that
$B_{R_0}(x_0)\subseteq\Omega$, that is $R_0\le
\dist(x_0,\partial\Omega)$\,. In the case $\Omega=\RR^d$ there is
obviously no restriction on $R_0$. Let $u_D$ be the solution to the
corresponding MDP, as defined in \eqref{Sub.Dirichlet.Problem}. It
has extinction time $T_m$.  By parabolic comparison, it is then
clear that $u_D\le u$ in $Q=(0,T')\times B_{R_0}(x_0)$, where
$T'=\min\{T,T_m\}$, hence we can easily extend the positivity
results for the MDP obtained in the previous section to any other
local weak solution, either in the form given by Theorem
\ref{Flux.Thm.Posit.2}, or in the Aronson-Caffarelli form
\eqref{AC.Estimates1}. This concludes the proof of Corollary
\ref{coroll1.2}, and as a consequence we get Theorem
\ref{Thm.Posit}.\qed

\subsection{\large Lower estimates independent of the extinction time for $m_c<m<1$.
Proof of Theorem \ref{GoodFDE.Posit}}

The proof  is divided in some short steps. Here, $m_c<m<1$.

\noindent$\bullet~$\textsc{Reduction.} Let $u_R(t,x)$ be the
solution to the homogeneous Dirichlet problem on the ball $B_R$,
corresponding to the initial datum $u_0\in\LL^1(B_R)$ and having
extinction time $T(R,u_0)$. The rescaled solution
\begin{equation}
{u}(t,x)=\frac{M}{R^d}\,\widehat{u}\left(\frac{t}{\tau},
\frac{x-x_0}{R} \right)\,,\qquad \tau= R^{2-d(1-m)}M^{1-m},\qquad
M=\int_{B_R}u_0(x)\,dx
\end{equation}
allows us to pass from a solution with mass $M$ defined in the ball
of radius $R$ centered at $x_0$ to a solution $\widehat{u}$ with
mass 1 in the ball of radius 1 centered at $0$. The extinction times
have to be rescaled accordingly,
$T(u)=R^{2-d(1-m)}M^{1-m}\overline{T}$, where
$\overline{T}={T}(\widehat u)$. Therefore, we will work with
rescaled problems and solutions.

\noindent$\bullet~$\textsc{Barenblatt solutions.} We now consider
the solution $\mathcal{B}$ of the Dirichlet problem posed on $B_1$,
and corresponding to the Dirac mass as initial data,
$\mathcal{B}(0,\cdot)=\delta_0$, and zero boundary data, that we
call Barenblatt solution. First we recall that by approximation with
$\LL^1$-data, or by comparison with the solutions of the Cauchy
problem, it is  easy to check that the smoothing effect applies and
\begin{equation}\label{Smooth.Baren.L1}
\mathcal{B}(t,x)\le c_{m,d}\,t^{-d\vartheta_1}, \qquad\mbox{for any
}(t,x)\in [0,+\infty)\times B_1\,.
\end{equation}
Moreover, it is known that this is the solution that extinguishes at
the later time among all nonnegative solutions with the same mass of
the initial data and same boundary data.  Such comparison is a
consequence of the concentration-comparison and symmetrization
arguments developed in detail in  \cite{VazBook, VazLN}. Thus, we
need to prove that the Barenblatt solution extinguishes in  finite
time $\overline{T}$. The proof is based on the fact that it is
bounded for $t\ge t_0>0$.

\noindent$\bullet~$\textsc{Solution by separation of variables}.
Consider now the solution
$\mathcal{U}_{r}(t,x)=S(x)(T_1-t)^{1/(1-m)}$ of the Dirichlet
problem on $B_r$, $r>1$. It corresponds to the initial datum
$\mathcal{U}_{r}(0,x)=S(x)T_1^{1/(1-m)}$, and extinguishes at a time
$T_1$, to be chosen later. Here, $S$ is the solution to the stationary
elliptic Dirichlet problem $\Delta S^m+{1/(1-m)}S=0$ on $B_{R_0}$,
and therefore it can be chosen radially symmetric, $S(x)=S(|x|)$. It
will also be nonincreasing in $r=|x|$. By standard regularity theory
$S(x)^m$ can be bounded from above an below by the distance to the
boundary. The parameter $T_1$, extinction time of $\mathcal{U}_{r}$
can be chosen at will. To fix it, we pick a time $t_0>0$ and define
the $T_1$ through the relation
\begin{equation}
S(1)(T_1-t_0)^{\frac{1}{1-m}}=c_{m,d}\,t_0^{-d\vartheta_1}.
\end{equation}

\noindent$\bullet~$\textsc{Comparing the two solutions.} We now
consider the homogeneous Dirichlet problem on $[t_0,T]\times B_{r}$,
and we compare the Barenblatt solution $\mathcal{B}$  with the
solution $\mathcal{U}_r$ constructed above in the cylinder
$Q=B_1(0)\times[t_0,T_1)$. At the initial time $t_0$ we know by
construction  that $\mathcal{B}(t_0,x)\le \mathcal{U}_r(t_0,x)$. The
comparison of the boundary data is immediate. By well-known
parabolic comparison results, this implies that $\mathcal{B}(t,x)\le
\mathcal{U}_r(t,x)$, on $[t_0,T_1]\times B_{r}$, and hence the
extinction times satisfy
\begin{equation}
\overline{T}\le
T_1=\left[\frac{c_{m,d}}{S(1)\,t_0^{d\vartheta_1}}\right]^{1-m}+t_0.
\end{equation}
We only need to choose a $t_0\in(0, T)$ to obtain an expression for
the upper bound of $\overline{T}$ that depends only on $m$ and $d$
(the reader may choose to optimize the expression for $T_1$ with
respect to $t_0$).

\noindent$\bullet~$\textsc{Conclusion. }  As a consequence of the
above upper bound, we know that any solution to the Dirichlet
problem on the unitary ball and with unitary initial mass extinguish
at a time $T\le \overline{T}\le \tau(m,d)$. Rescaling back to the
original variables we have proved that any solution $u_R$ to the
Dirichlet problem on the ball $B_R$ and with initial mass
$M=\int_{B_R}u_0\dx$ extinguish at a time $T_R$ that can be bounded
above by
\[
T(R,u_0)\le \tau_{m,d}R^{2-d(1-m)}\|u_0\|_{\LL^1(B_R)}^{1-m}.
\]
The lower bounds come from the fact that $t_*\le T$ and is given by
\eqref{form.t_*}.

\noindent$\bullet~$ \textsc{The lower Harnack inequality}.
Inequality  \eqref{GoodFDE.Low.Harnack} follows by plugging the
upper bound \eqref{GoodFDE.sharp.ext.est} into the lower bound
\eqref{Flux.Posit.Sub.Dirichlet.2}.\qed

\medskip

The reader should notice that the properties that we have used are
typical of the good fast diffusion range, $m_c<m<1$, and cannot be
extended to the very fast diffusion range, $m<m_c$.


\subsection{\large Lower estimates independent of the extinction for $0<m<m_c$.}
\label{sect-loindT}

The presence of the minimal extinction time $T_m=T$ in the formula
for the lower Harnack inequality responds to an essential characteristic of the
problem. Actually, lower estimates in terms of only $\LL^1$ norms
cannot be true for $m\le m_c$ as we have shown at the beginning of this section:
there is no positive lower bound at a time $t_0>0$ and a point $x_0$
that depends only on $t_0, R$ and the mass of $u_0$ in $B_R(x_0)$.
Similar examples can be constructed if $u_0\in \LL^p_{loc}(\RR^d)$
with $p <p_c$, and we can be found in \cite{VazLN}, Chapters 5 and~
7.

Fortunately, controlling the local (or global) $\LL^{p}$ norm gives
a control on the MET $T$, and in this way we get valid lower
estimate without $T$, as we explain next.

\subsubsection{Estimates in terms of the $\LL^{p_c}$
norm}\label{sec.upper.fet.1}

In  \cite{BCr-cont}  B\'enilan and Crandall prove that for any $0\le
s\le t\le T$, and for any $m<m_c$
\begin{equation}\label{pc.Norm.Decay}
\|u(t)\|_{p_c}^{1-m}\,\le\,\|u(s)\|_{p_c}^{1-m}-\,\mathcal{K}_{p_c}\,(t-s)\;,
\qquad\mbox{with}\qquad\mathcal{K}_{p_c}=\frac{8\,[d(1-m)-2]\,\mathcal{S}_2^2}{(d-2)^2(1-m)}\;,
\end{equation}
where $\mathcal{S}_2$ in the constant of the Sobolev inequality
\begin{equation}\label{Sob.2}
\|f\|_{2^*}\le\,\mathcal{S}_2\,\|\nabla f\|_2
\end{equation}
and the above estimate holds for any solution with initial datum
$u_0\in\LL^{p_c}$. We also stress on the fact that the constant
$\mathcal{K}_{p_c}$ is universal in the sense that it only depends
on $m$ and $d$. As a consequence of \eqref{pc.Norm.Decay}\,, we
have the following \textsl{universal upper bound for the
extinction time}
\begin{equation}\label{pc.Norm.Ext.Time}
T(u_0)\,\le\,\mathcal{K}_{p_c}^{-1}\,\|u_0\|_{p_c}^{1-m}\;.
\end{equation}
\textit{We remark that while for lower bounds on FET we only need local
information on the initial datum, upper estimates for the FET
require  global information. Fortunately, in the minimal problem
that we are considering, global and local are equivalent since
$u_0(x)=0$ for $|x-x_0|\ge R$.}

\medskip

\noindent {\sc Proof.~} We sketch here the proof for the reader's
convenience. It is well known that the time derivative of the global
$\LL^p$ norm of the solution $u(t)$ of the MDP problem under
consideration is given by
\begin{equation}\label{Lp.Global.Time.Deriv}\begin{split}
\frac{\rd}{\dt}\|u(t)\|_p^p
&=-\frac{4p(p-1)}{(p+m-1)^2}\int\left|\nabla u^{\frac{p+m-1}{2}}\right|^2\dx\\
&\le -\frac{4p(p-1)\,\mathcal{S}_2^2}{(p+m-1)^2}\left[\int
u^{\frac{(p+m-1)2^*}{2}}\dx\right]^{\frac{2}{2^*}}
= -\frac{4p(p-1)\,\mathcal{S}_2^2}{(p+m-1)^2}\,\|u\|_{\frac{(p+m-1)2^*}{2}}^{p+m-1}\;,\\
\end{split}
\end{equation}
where in the last step we used the Sobolev inequality \eqref{Sob.2}
applied to the function $f=u^{(p+m-1)/2}$, where $2^*=2d/(d-2)$
and $\mathcal{S}_2$ is the Sobolev constant.

Notice that if $m>m_c$, then $p_c<1$, so that the global
$\LL^p$-norm increases, and this originates the so called Backward
Effect, see e.g. \cite{VazLN}. This explains our assumption
$m<m_c$\,. Moreover
\[
p_c+m-1=p_c\left(1-\frac{2}{d}\right), \qquad
\frac{(p_c+m-1)2^*}{2}=p_c\,,\qquad
\frac{4p_c(p_c-1)\,\mathcal{S}_2^2}{(p_c+m-1)^2}
=\frac{8\,[d(1-m)-2]\,\mathcal{S}_2^2}{(d-2)^2(1-m)}>0,
\]
so that \eqref{Lp.Global.Time.Deriv} becomes
\[
\frac{\rd}{\dt}\|u(t)\|_{p_c}^{p_c}\,\le\,-\frac{8\,[d(1-m)-2]\,\mathcal{S}_2^2}{(d-2)^2(1-m)}\|u(t)\|_{p_c}^{p_c\,\left(1-\frac{2}{d}\right)}
\]
integrating the differential inequality gives the bound \eqref{Lp.Global.Time.Deriv}
for any $0\le s \le t$\,. Letting $s=0$ and $t=T(u_0)$ in \eqref{Lp.Global.Time.Deriv} finally gives
\eqref{pc.Norm.Ext.Time}.\qed

\medskip

\noindent {\bf Application to Theorem \ref{Thm.Posit.v2}.} Using
this bound, we can now formulate the second version of the lower
Harnack estimate, reflected in the theorem. The proof is immediate
when $u_0$ is as in the minimal problem, since in that case local
and global norm is the same. Comparison as done in Subsection
\ref{sect.pls}, allows to pass to the general solutions. Notice
that in this way we use a local $\LL^{p_c}$ norm, not the global
one!

\medskip

\noindent {\bf Remark.} When we have not only the Sobolev
inequality, but also the Poincar\'e, we can prove similar estimates
for any $m\in(0,1)$\,.  This happens for instance for problems posed
in bounded domains, or for the minimal Dirichlet problem.

\subsubsection{Estimates in terms of other $\LL^{p}$
norms}\label{sec.upper.fet.2}

\begin{prop}\label{ext.dirich}
Let $m<1$\,, $\alpha\ge 1$\,, $R>0$ and let $u$ be the solution to
the Dirichlet problem
\begin{equation*}
\left\{\begin{array}{lll}
 u_t=\frac{1}{m}\Delta (u^m) & ~ {\rm in}~ (0,T)\times B_{\alpha R}\\[2mm]
 u(0,x)=u_0(x) & ~{\rm in}~ B_{\alpha R},~~~{\rm and}~~~ \supp(u_0)\subseteq B_R\\[2mm]
 u(t,x)=0 & ~{\rm for}~  t >0 ~{\rm and}~ x\in\partial B_{\alpha R}\\
\end{array}\right.
\end{equation*}
with $u_0\in\LL^p\big(B_{\alpha R}\big)$, with
$p>\max\{p_c,1\}=\max\left\{d(1-m)/2,1\right\}$\,. Then the following estimate
\begin{equation}\label{Global.Lp.Estimates}
\big\|u(t)\big\|_p^{1-m}-\big\|u(s)\big\|_p^{1-m}
 \le -\mathcal{K}_p\;(t-s)\;.
\end{equation}
hold for any $0\le s\le t$\,, where
\[
\mathcal{K}_p=4 \frac{(1-m)(p-1)}{\big[p+m-1\big]^2}
 \left[\mathcal{P}\,\alpha R\,\right]^{-2\left(1-\frac{p_c}{p}\right)}\;
 \mathcal{S}_2^{-\frac{2p_c}{p}}>0
\]
and where $\mathcal{S}_2$ is the Sobolev constant of $\RR^d$ and
$\mathcal{P}$ is the Poincar\'e constant on the unit ball\,.
\end{prop}

\noindent {\sl Proof.~} First consider, for any $f\in W_0^{1,2}\big(B_{\alpha R}\big)$\,,
the Sobolev and Poincar\'e inequalities:
\[
\|f\|_{2^*}\;\le\;\mathcal{S}_2\;\|\nabla\!f\|_2\,,
    \qquad\mbox{and}\qquad
    \|f\|_{2} \;\le\;\mathcal{P}\;\alpha\,R\;\|\nabla\!f\|_2\;,
\]
where $2^*=2d/(d-2)$\,, and where the constants $\mathcal{S}_2$, the
optimal Sobolev constant on $\RR^d$,  and $\mathcal{P}$, the
Poincar\'e constant on the unit ball, only depend on the dimension
$d$\,. By combining them through the H\"older inequality, we then
get for any $q\in(2,2^*)$
\[
\big\|f\big\|_{q}
 \le \big\|f\big\|_{2}^{1-\vartheta}\big\|f\big\|_{2^*}^\vartheta
 \le \left[\mathcal{P}\,\alpha R\,\right]^{1-\vartheta}
    \;\mathcal{S}_2^\vartheta\;\big\|\nabla\!f\big\|_2\,.
\]
Let now
\[
f=u^{\frac{p+m-1}{2}}\;,\qquad q:=\frac{2p}{p+m-1}
\qquad\mbox{and}\qquad\vartheta=\frac{d(1-m)}{2p}=\frac{p_c}{p}\,.
\]
We remark that $q<2^*$ if and only if $p>p_c$, while $q>2$ if and
only if $q<\infty$\,. We obtain then
\begin{equation}\label{Sob+Poinc}
\left\| u \right\|_p^{p\left[1-\frac{1-m}{p}\right]}
 \le \left[\mathcal{P}\,\alpha R\,\right]^{2\left(1-\frac{p_c}{p}\right)}\;
  \mathcal{S}_2^{\frac{2p_c}{p}}\;
  \left\|\nabla\!u^{\frac{p+m-1}{2}} \right\|_2^2
 :=\mathcal{K}_0\left\|\nabla\!u^{\frac{p+m-1}{2}} \right\|_2^2
\end{equation}
The derivative of the global $\LL^p$-norm then satisfies
\begin{equation}
\frac{\rd}{\dt}\left\| u(t) \right\|_p^p
 =-\frac{4p(p-1)}{\big[p+m-1\big]^2}\left\|\nabla\!u^{\frac{p+m-1}{2}} \right\|_2^2
 \le-\frac{4p(p-1)\mathcal{K}_0^{-1}}{\big[p+m-1\big]^2}\big\|u\big\|_p^{p\left[1-\frac{1-m}{p}\right]}
\end{equation}
where in the last step we used \eqref{Sob+Poinc}. Integrating the differential inequality over $[s,t]\subseteq
[0,T]$, gives
\[
\big\|u(t)\big\|_p^{1-m}-\big\|u(s)\big\|_p^{1-m}
 \le -\frac{4(1-m)(p-1)}{\big[p+m-1\big]^2}\left[\mathcal{P}\,\alpha R\,\right]^{-2\left(1-\frac{p_c}{p}\right)}\;
  \mathcal{S}_2^{-\frac{2p_c}{p}}\;(t-s)\;.\mbox{~\qed}
\]

\medskip

\noindent \textbf{Upper Bounds on the Extinction Time.} The above
estimates \eqref{Global.Lp.Estimates}, prove that any solution of
the Dirichlet problem extinguish in finite time, and this is not
surprising, but they also provide an \textit{Upper Bound for the
extinction time} $T$\,, indeed letting $s=0$ and $t=T$, we obtain
\begin{equation*}
T \le \mathcal{K}_p^{-1}\big\|u_0\big\|_p^{1-m}
  = \frac{\big[p+m-1\big]^2}{4(1-m)(p-1)}
 \left[\mathcal{P}\,\alpha R\,\right]^{2\left(1-\frac{p_c}{p}\right)}\;
 \mathcal{S}_2^{\frac{2p_c}{p}}\big\|u_0\big\|_p^{1-m}
\end{equation*}
Notice that in the limit $p\to p_c$ we recover the previous result
\eqref{pc.Norm.Decay}. Summing up, the above result proves that a
\textit{ global Sobolev and Poincar\'e inequality provides that
the solution extinguish in finite time $T$ and an gives a
quantitative upper bound for $T$\,}.

\medskip

\noindent {\bf Remarks.} These results can be extended to
different domains or manifolds in a straightforward way, the only
important thing is to have global Sobolev and Poincar\'e
inequalities, as already studied by the authors in \cite{BGV-JEE},
in the case of Riemannian manifolds with nonpositive curvature

\medskip

Using this bound, we can now formulate a  version of the lower
Harnack estimate similar to Theorem \ref{Thm.Posit.v2}. We leave the
easy details to the reader.

\section{\large  Part II\,. Local upper bounds}

In  the second part of this work we turn our attention to the
question of upper estimates for solutions with data in some
$\LL^p_{loc}$,  $p\ge 1$,  and obtain quantitative forms of the
bounds that are sharp in various respects. The range of
application is all $m<1$, even $m\le 0$. We assume moreover that
$d\ge 3$, which is the interesting case also for the lower
estimates, in order to avoid technical complications which break
the flow of the proofs and results, but we remark that the
qualitative fact, the existence of local upper bounds, is also
true for $d=1,2$.

\smallskip

 As a preliminary for  the main result, we
devote Section \ref{sect.lp} to establish the conservation of the
local $\LL^p$ integrability of the solutions and the control of
the evolution of the local $\LL^p$ norm for suitable $p\ge 1$. Let
$u=u(t,x)$ be a nonnegative weak solution of the FDE for $m<1$
defined in  a space-time cylinder $Q=(0,T]\times B_{R_0}$ for some
$R_0,T>0$. This is the form of the estimate we get:
\begin{equation*}
\left[\int_{B_R(x_0)}|u(t,x)|^p \dx \right]^{(1-m)/p}\le
\left[\int_{B_{R_0}(x_0)}|u(s,x)|^p\dx\right]^{(1-m)/p} +
K\;(t-s),
\end{equation*}
for any $R_0>R$ and $0\le s\le t< T$. It is valid for all $m<1$ if
$p\ge 1$, $p>1-m$. The dependence of $K$ on $R$ and $R_0$ is
explicitly given in Theorem~\ref{Lem.Local.Lp.norms} below. The
estimate extends Herrero-Pierre's well-known estimate to $p>1$ and
is valid for $m\le 0$.

\smallskip

 The main result of this part is   the local
upper bound  that applies for the same type of solution and
initial data, under different restrictions on $p$. Here is the
precise formulation.

\begin{thm} \label{thm.upper} Let $p\ge 1$ if $m>m_c$
or  $p>p_c$ if $m\le m_c$. Let $u$ be a local weak solution to the FDE in the cylinder $(0,T)\times\Omega\subseteq (0,+\infty)\times\RR^d$.
Then there are positive constants $\mathcal{C}_{1}$,
$\mathcal{C}_{2}$ such that we have
\begin{equation}\label{upper.est}
\begin{split}
u(t,x_0)&\le\frac{\mathcal{C}_{1}}{t^{d\vartheta_p}}\,\left[
\int_{B_{R_0}(x_0)}|u_0(x)|^{p}\dx
\right]^{2\vartheta_p}+\mathcal{C}_{2}\left[\frac{t}{R_0^2}\right]^{\frac{1}{1-m}}.
\end{split}
\end{equation}
where $R_0\le\dist(x_0,\partial\Omega)$ and
the constants $\mathcal{C}_{i}$ depend on $m,d$
and $p$.
\end{thm}

 We recall that $\vartheta_p=1/(2p-d(1-m))=1/2(p-p_c)$.
Note that the constants $\mathcal{C}_i$ do not depend on the radii, but only
on $m,\,d$ and $p$. An explicit formula for them is given at the end of
the proof, but we point out that such values need not be optimal.
 The result is
proved in Sections \ref{sect3} and \ref{sect4}. A similar smoothing
effect result has been proved for the first time by Herrero and
Pierre in \cite{HP} in the good fast diffusion range $m_c<m<1$ using
$p=1$, but it is new in the range $m\le m_c$ where HP's result
cannot hold in view of solutions like \eqref{VSS}. HP's technique
relies on stronger differential estimates that do not hold in the
subcritical fast diffusion case or on the local setting; our
impression is that their techniques can not be adapted to the very
fast diffusion range. Related estimates for $p>1$ are due to
DiBenedetto and Kwong, \cite{DK}\,, and Daskalopoulos and Kenig,
\cite{DasKe}, but as far as we know no results cover the very fast
diffusion range.  Finally, note  that the smoothing effect
$\LL^p_{loc}$ into $\LL^\infty_{loc}$ is false for exponents $p<p_c$ as
has been demonstrated in \cite{VazLN}. In fact, that monograph
studies the existence of the so-called backward smoothing effects
that go from $\LL^p(\RR^d)$ into $\LL^1(\RR^d)$ for $p<p_c$.

The local bound in (\ref{upper.est}) is expressed as the sum of two
independent terms, one due to initial data, the other one due to
effects near the boundary. The estimate is optimal in the following
senses:

\noindent (i) The first term responds to the influence of the
initial data and has the exact form that has been demonstrated to
be exact for solutions that are defined in the whole space and
have initial data in $\LL^p(\RR^d)$, see \cite[Chapters 3,5]{VazLN}.
By exact we mean that the integral is the same (but extended to
the whole space) and the exponents are the same, only the constant
$\mathcal{C}_{1}$ may differ.  We can then recover the global
smoothing effect on $\RR^d$, just by letting $R\to\infty$ so that
the second term disappears; as mentioned above the constant
$\mathcal{C}_1$ is not the optimal one: the best constant for the
global smoothing effect on $\RR^d$ has been calculated by one of
the authors in \cite{VazLN}\,.

\noindent (ii) The last term accounts for the influence of the
boundary data and is special to fast diffusion in the sense that it
does include any information on the precise boundary data, thus
allowing for the so-called {\sl large solutions} that take on the
value $u=+\infty$ on $\partial B_R$. The term has the exact form
prescribed by the explicit singular solutions \eqref{VSS}. This last
term has the meaning of an absolute bound for all solutions with
zero or bounded initial data; thus it can also be interpreted as a
{\sl universal bound} for the influence of any boundary effects.
 In applications it is interpreted as an absolute
damping of all external influences.

\subsection{\large Evolution  of Local $\LL^p$-norms}
\label{sect.lp}

A basic question in the existence theory is obtaining a priori
estimates of the solutions in terms of the data measured in some
appropriate norm. The peculiar feature of the FDE is the local
nature of the estimates. A fundamental result in this direction is
the local $\LL^1_{\rm loc}$-$\LL^1_{\rm loc}$ estimate due to Herrero
and Pierre (which is valid for $m>0$):

\begin{lem}\label{Lem.Local.L1.norms}
Let $u,v\, \in\,C\left([0,+\infty)\,;\,\LL^1_{\rm
loc}(\RR^d)\right)$ be weak solutions of
\begin{equation*}
\partial_t u =\Delta(u^m/m), \quad 0<m<1.
\end{equation*}
Let $R>0$, $R_0=\lambda R$ with $\lambda>1$, and $x_0\in \RR^d$ be
such that $B_1=B_{R_0}(x_0)\subset \RR^d$.  Let moreover $v\le u $
a.e. Then, the following inequalities hold true:
\begin{equation}\label{Lem.Local.Lp.p=1}
\left[\int_{B_R}\big[u(t,x)-v(t,x)\big] \dx \right]^{1-m}\le
\left[\int_{B_1}\big[u(s,x)-v(s,x)\big]\dx\right]^{1-m} +
K_{R,R_0, 1}\;|t-s|,
\end{equation}
for any $t,s\ge 0$, where
\begin{equation}
K_{R,R_0, 1}= \frac{c_1}{\big(R_0-R\big)^2}\vol\left(B_{
R_0}\setminus B_R\right)^{(1-m)}>0
\end{equation}
and the constant $c_1>0$ depends only on $m$, $d$.
\end{lem}

This result  was proven in  Prop. 3.1 of \cite{HP} and has been
generalized to the case of fast diffusion on a Riemannian manifold
by the authors in \cite{BGV-JEE}. Our goal here is to extend such
a result into and $\LL^p_{\rm loc}$-$\LL^p_{\rm loc}$ estimate for
suitable $p>1$. This estimate has two merits: first, it is valid
for all $-\infty<m<1$; second, it is needed for some values
$p>p_c$ for the proof of boundedness estimates.

\bigskip

\begin{thm}\label{Lem.Local.Lp.norms}
Let $u\, \in\,C\left((0,T)\,;\,\LL^1_{\rm loc}(\Omega)\right)$ be
a nonnegative weak solution of
\begin{equation}\label{Def.Distribution.Solution.}
\partial_t u =\Delta(u^m/m),
\end{equation}
and assume that  $u(t,\cdot)\in\LL^p_{\rm loc}(\Omega)$ for some
$p\ge 1$, $p>1-m$, and for all $0< t< T$. Here, $\Omega$ is a
domain in $\RR^d$ that contains the ball $B_1=B_{ R_0}(x_0)$.
Then, the following inequality holds true:
\begin{equation}
\left[\int_{B_R(x_0)}|u(t,x)|^p \dx \right]^{(1-m)/p}\le
\left[\int_{B_{R_0}(x_0)}|u(s,x)|^p\dx\right]^{(1-m)/p} +
K_{R,R_0, p}\;(t-s),
\end{equation}
for any $0\le s\le t< T$, where
\begin{equation}
K_{R,R_0, p}= \frac{p\,c_{m,d}}{\big(R_0-R\big)^2}\,\vol\left(B_{
R_0}\setminus B_R\right)^{(1-m)/p}>0,
\end{equation}
and the constant $c_{m,d}>0$ depends only on $m$, $d$.
\end{thm}

\noindent {\bf Remarks.} (i) The result implies  for those values of
$p$ that  whenever $u(s,\cdot)\in\LL^p_{\rm loc}(\Omega)$ for some
$s>0$, then $u(t,\cdot)\in\LL^p_{\rm loc}(\Omega)$ for all $t>s$.
Note that the dependence of the local $\LL^p$ norm is again expressed
as the sum of two independent terms, one due to the initial data,
the other one due to effects near the boundary.

\noindent (ii)  Note that  the times $t$ and $s$ must be ordered
in this result, a condition that  is not required in Lemma
\ref{Lem.Local.L1.norms}.

\noindent (iii) The last term may blow up as we approach the
boundary of $\Omega$ (where no information on the data is used).
Indeed, the constant can be written in the form
$$
K_{R,R_0,p}=pc'_{m,d}R_0^{2(p-p_c)/p}F(R/R_0), \quad F(s)=
\frac{(1-s^d)^{(1-m)/p}}{(1-s)^2}.
$$
Now, if $x_0\in
\Omega$ we may take $R_0=d(x_0,\partial \Omega)$ and  $R=R_0(1-\ve)$.
In that case the constant in the last term behaves as $\ve\to 0$
in the form
$$
K_{R,R_0,
p}\sim R_0^{2(p-p_c)/p}\ve^{-\beta}, \quad \beta=2-(1-m)/p=(2p+m-1)/p.
$$

\noindent (iv) The constant blows up in the limit $m\to 1^-$\,,
and this is perfectly coherent, since a similar estimate is false
for the Heat Equation.

Moreover, the constant $K_{R,R_0, p}$, blows up when $p\to\infty$,
thus it does not provide $\LL^\infty$ local stability, while it
provides local $\LL^p$ stability, for $p>p_c$.

\medskip

\noindent {\sc Proof of Theorem \ref{Lem.Local.Lp.norms}.~} (i)
Let $u \ge 0$ and take a test function   $\psi\in
C_c^\infty(\Omega)$  and $\psi\ge 0$. We can compute
\begin{equation}\label{dt.Psi}
\begin{split}
\frac{\rd}{\dt}\int_{{\Omega}} \psi\, u^p \dx &=
p\int_{{\Omega}}\psi\, u^{p-1}\partial_t u
\dx=-p\int_{{\Omega}}
\nabla\left(\psi\, u^{p-1}\right)\cdot\nabla\left(\frac{u^m}{m}\right) \dx\\
&=-p\left[\int_{_{\Omega}}\nabla\psi\cdot(u^{p+m-2}\nabla u)\dx+(p-1)\int_{{\Omega}}\psi\, u^{p+m-3}|\nabla u|^2 \dx\right]\\
&=-p\left[\frac{1}{p+m-1}\int_{{\Omega}}\nabla\psi\cdot\nabla(u^{p+m-1})\dx
+\frac{4(p-1)}{(p+m-1)^2}\int_{{\Omega}}\psi\,\left|\nabla u^{\frac{p+m-1}{2}}\right|^2 \dx\right]\\
&=\frac{p}{p+m-1}\int_{{\Omega}}(\Delta\psi)\,u^{p+m-1}\dx
-\frac{4p\,(p-1)}{(p+m-1)^2}\int_{{\Omega}}\psi\,\left|\nabla u^{\frac{p+m-1}{2}}\right|^2 \dx\\
&\le\frac{p}{p+m-1}\int_{{\Omega}}\big|\,\Delta\psi\,\big|\,u^{p+m-1}\dx.
\end{split}
\end{equation}
 This computation
holds true for any $p\ge1$,  and any $m\in\RR$, when one replaces,
in the limit $m\to 0$, the quantity $(u^m)/m$ with  $\log u$; we
also have to replace $u^{p+m-1}/(p+m-1)$ by $\log(u)$ if
$p+m-1=0$. Of course, when $p+m-1\le 0$ the last term may be
infinite, since it contains $u^{p+m-1}$ so we make the assumption  $p>1-m$.

\noindent (ii) Under such  assumptions, inequality \eqref{dt.Psi}
implies that for any solution $u\ge 0$ we have

\begin{equation}\label{Lem.HP.2}
\begin{split}
\frac{\rd}{\dt}\int_{{\Omega}}\psi(x)  u(t,x)^p\dx &\leq
\frac{p}{p+m-1}\int_{\Omega} |\Delta(\psi(x))|\; |u(t,x)|^{p+m-1}\dx\\
&\le \frac{p}{|p+m-1|}
\left[\int_{{\Omega}}|\Delta(\psi(x))|^{\frac{p}{(1-m)}}\;
\psi(x)^{1-p/(1-m)} \dx \right]^{\frac{1-m}{p}}
\left[\int_{{\Omega}}\psi(x)\, u(t,x)^p\dx \right]^{\frac{p+m-1}{p}}\\
&=C(\psi)\left[\int_{{\Omega}}\psi(x)\, u(t,x)^p\dx
\right]^{1-\frac{1-m}{p}}
\end{split}
\end{equation}
where in the second step we have used H\"older's inequality with
conjugate exponents $p/(1-m)$ and $p/(p+m-1)$, and where
\begin{equation}\label{C.Psi}
C(\psi)=\frac{p}{p+m-1}
\left[\int_{\Omega}|\Delta(\psi(x))|^{\frac{p}{(1-m)}}\psi^{1-\frac{p}{(1-m)}}
\dx \right]^{\frac{1-m}{p}}.
\end{equation}
 We will check below that this quantity can be made
finite by a proper choice of $\psi$. Assuming this for the moment,
formula \eqref{Lem.HP.2} can be expressed as  a differential
inequality of the form
\[
y^\prime(\tau)\le C y^{1-\varepsilon}(\tau)
\]
where $y(\tau)=\int_{_{\Omega}}\psi(x)|u(\tau,x)|^p\dx$,
$C=C(\psi)$ and $\varepsilon=(1-m)/p\in (0,1).$ Integrating such
differential inequality over $(s,t)$ lead to
\[
y^\varepsilon(t)-y^\varepsilon(s)\le C\,\varepsilon\, (t-s)
\]
that is
\[
\left[\int_{_{\Omega}}\psi(x)|u(t,x)|^p\dx\right]^{(1-m)/p}\le
\left[\int_{\Omega}\psi(x)|u(s,x)|^p\dx\right]^{(1-m)/p}+\frac{(1-m)}{p}\,C(\psi)(t-s)
\]
for any $0\le s\le t$. This will immediately imply the statement,
once we prove the bounds
\begin{equation}\label{Bounds.C(psi)}
\frac{(1-m)}{p}\,C(\psi)= K_{R,\lambda,p}<+\infty.
\end{equation}

\noindent (iii) We only have to verify the form of the last bound.
To this end we consider a function $\psi=\varphi^b\in
C_c^{\infty}(M)$, with
\begin{equation}\label{Lem.HP.7}
0\le \varphi\le 1\,,\quad \varphi\equiv 1 ~\mbox{in~} B_R\,,\quad
\varphi\equiv 0 ~\mbox{outside~} B_{\lambda R}
\end{equation}
with $\lambda=R_0/R>1$. Moreover, we will assume that $\varphi$ is
radially symmetric and $
\varphi(x)=\overline{\varphi}\left(|x|/R\right)$\,,
where $\overline{\varphi}:\RR\to\RR$ is a $C_c^{\infty}(\RR)$
function such that:
\[
0\le \overline{\varphi}(s)\le 1\,,\quad
\overline{\varphi}(s)\equiv 1\,, ~\mbox{for~} 0\le s\le 1\,,\quad
\overline{\varphi}\equiv 0\,, ~\mbox{for~} s\ge\lambda
\]
where $\lambda>1$ and $|x|$ is the distance from a fixed point. We
then have
\begin{equation}\label{Bounds.C(psi).2}\begin{split}
\left|\Delta\left(\psi(x)\right)\right|^{\frac{p}{1-m}}&\psi(x)^{1-\frac{p}{1-m}}=
\varphi(x)^{b\left[1-\frac{p}{1-m}\right]}\left|b(b-1)\,\varphi^{b-2}\,
\left|\nabla\,\varphi\,\right|^2+b\,\varphi^{b-1}\,\Delta\varphi\right|^{\frac{p}{1-m}}\\
\le &
[b(b-1)]^{\frac{p}{1-m}}\varphi^{b\left[1-\frac{p}{1-m}\right]+\frac{(b-2)p}{1-m}}
\left|\,\left|\nabla\,\varphi\,\right|^2+\left|\Delta\varphi\right|\right|^{\frac{p}{1-m}}
\end{split}
\end{equation}
the last inequality follows from the fact that we are considering
a radial function $0\le \varphi(x)=\overline{\varphi}(|x|/R)\le
1$, with $b>\tfrac{2p}{1-m}$. Then, we compute
\[\begin{split}
|\nabla\varphi(x)|^2 &=R^{-2}|\overline{\varphi}'(|x|/R)|^2\,|\nabla |x||^2\le R^{-2}
\,|\overline{\varphi}'(|x|/R)|^2
\le c'^{2}_\lambda\,R^{-2}\\
|\Delta\varphi(x)|&=\left|R^{-2}\overline{\varphi}''(|x|/R)\,|\nabla(|x|)|^2+
R^{-1}\overline{\varphi}'(|x|/R)\Delta|x|\right|\\
&\le\frac{1}{R}\,\left[\,\frac{|\overline{\varphi}''(|x|/R)|}{R}+
|\overline{\varphi}'(|x|/R)|\,\frac{d-1}{|x|}\,\right]
\le\frac{(d-1)c''_\lambda}{R^2},
\end{split}
\]
where in the last step we used the fact that $\Delta\varphi$ is
supported in $A_{R\,,\lambda}=B_{\lambda R}\setminus B_R$ and that
the smooth function $\overline{\varphi}$ has bounded derivatives
in $A_{R,\lambda}$
\begin{equation}\label{Bounds.C(psi).3}
\big|\overline{\varphi}(|x|/R)\big|\le
\frac{c'_0}{\lambda-1}=c'_\lambda\;,
\qquad\big|\,\overline{\varphi}'(|x|/R)\big|
+\big|\overline{\varphi}''(|x|/R)\big|\le
\frac{c''_0}{(\lambda-1)^2}=c''_\lambda\;,
\end{equation}
we just remark that this last estimate depend on an explicit choice
of the test function $\overline{\varphi}$.\\
Inequality \eqref{Bounds.C(psi).2} together with
\eqref{Bounds.C(psi).3} gives
\begin{equation*}\begin{split}
\left|\Delta\left(\psi(x)\right)\right|^{\frac{p}{1-m}}\psi(x)^{1-\frac{p}{1-m}}
&\le\,
[b(b-1)]^{\frac{p}{1-m}}\varphi^{b\left[1-\frac{p}{1-m}\right]+\frac{(b-2)p}{1-m}}
\left|\,\left|\nabla\,\varphi\,\right|^2+\left|\Delta\varphi\right|\right|^{\frac{p}{1-m}}\\
&\le[b(b-1)]^{\frac{p}{1-m}}\left[\frac{c'^2_0+(d-1)c''_0}{\big[(\lambda-1)R\big]^2}\right]^{\frac{p}{1-m}}
:=\frac{c'^{\frac{p}{1-m}}_p}{\big[(\lambda-1)R\big]^{\frac{2p}{1-m}}}
\end{split}
\end{equation*}
if $b>\frac{2p}{1-m}$, $c'_p=b(b-1)(c'^2_0+(d-1)c''_0)$. An
integration over $A_{R,\lambda}$ gives:
\begin{equation*}\begin{split}
\frac{(1-m)}{p}\,C(\psi)& =\frac{1-m}{p+m-1}
\left[\int_{A_{R,\lambda}}|\Delta(\psi(x))|^{\frac{p}{(1-m)}}\psi^{1-\frac{p}{(1-m)}} \dx \right]^{\frac{1-m}{p}}\\
&\le\frac{1-m}{p+m-1}\frac{c'_p}{\big[(\lambda-1)R\big]^2}\vol(A_{R,\lambda})^{\frac{1-m}{p}}
=\frac{c_p}{\big[(\lambda-1)R\big]^2}\vol(A_{R,\lambda})^{\frac{1-m}{p}}\,:=K_{R,\lambda,p}<+\infty
\end{split}
\end{equation*}
where $c_p=b(b-1)(c'^2_0+(d-1)c''_0)(1-m)/(p+m-1)$, and
$b>2p/(1-m)$, we can choose $b=3p/(1-m)$ to get
\[
c_p\le c_{m,d}\,p
\]
where $c_{m,d}$ is independent of $p$. The proof is thus
complete.\qed

\subsection{\large Smoothing effect in terms of space-time integrals}
\label{sect3}

 In this section we are going to prove  a first version
of the Local Smoothing Effect for the FDE. More precisely, we
prove that $\LL^p_{\rm loc}$ regularity in space-time implies
$\LL^\infty_{\rm loc}$ estimates, even when $m<m_c$, on the
condition that $p$ must be large enough. The estimates are local,
both in space and in time, but uniform on balls and the dependence is quantitative.
We consider a nonnegative weak solution of the FDE for $m<1$ defined
in  a space-time cylinder $Q=(0,T]\times B_R$ for some $R,T>0$.

 Throughout this section $T$ will not denote extinction time.

\begin{thm}\label{Thm.1.Integrated} Let $u$ and $m$ be as above, and let $p\ge 1$ if $m>m_c$
or  $p>p_c$ if $m\le m_c$.  For any two finite cylinders
$Q_1\subset Q_0$, $Q_i=(T_i,T]\times B_{R_i}$, with $0<R_1<R_0$,
and $0\le T_0<T_1<T$, we have
\begin{equation}
\sup_{Q_1}|u|\,\le\,\mathcal{C}_{\rm loc}\,
        \left[\frac{1}{(R_0-R_1)^2}+\frac{1}{T_1-T_0}\right]^{\frac{d+2}{2p+d(m-1)}}
        \;\left[\iint_{Q_{0}}u^{p}\dx \dt+\vol(Q_0)\right]^{\frac{2}{2p+d(m-1)}}\,.
\end{equation}
Moreover, the constant $\mathcal{C}_{\rm loc}$ depends only on
$m,d,p$.
\end{thm}

The proof presented here uses Moser's iteration process, and
borrows  some ideas  of \cite{DasKe} and
\cite{DiBenedettourbVesp}.  We will consider nested space-time
cylinders, in order to obtain the first estimates needed to prove
local SE. The proof will consist of the combination of several
partial results, which maybe of independent interest, and will be
split into several steps.  Note that by scaling the proof of this
kind of result need only to be done for a unit cylinder $Q_0$
where $R_0=1$ and $T_1-T_0=1$, and this is the case that will be
needed in the sequel.

\subsubsection*{ Step 1. Space-Time Energy Inequality}
Now we consider a solution $u$ defined in a parabolic cylinder
$Q=(T_0,T]\times B_R$ for some $R > R_1>0$, $T>0$ and consider
another parabolic cylinder $Q_1=(T_1,T]\times B_{R_1}$\,,
contained in $Q$, since we also let $T_0<T_1<T$. Then

\begin{lem}\label{Lem.Step.1} Under these assumptions, for every $m<1$ and $p>\max\{1,1-m\},$ we have
\begin{equation}\label{Step.1.final}
\begin{split}
\int_{B_{R_1}}u^p(T,x) \dx &+\iint_{Q_1}
\,\left|\nabla
u^{\frac{p+m-1}{2}}\right|^2\,\dx \dt
\le \mathcal{C}(m,p)\,\left[
\frac{1}{(R-R_1)^2}+\frac{1}{T_1-T_0}\right]
\left[\iint_{Q}\,\left(u^{p+m-1}+\,u^{p}\,\right)\dx \dt\right].
\end{split}
\end{equation}
The result also holds  when $u$ is a sub-solution, i.e.,  $u_t\le
\Delta u^m$\,.
\end{lem}
\noindent {\sl Proof.~} (i) We multiply the equation $\partial_t u
=\frac{1}{m}\Delta u^m$ by $\psi^2\,u^{p-1}$, with $p>1$ to be
chosen later in a suitable way, we take $\psi=\psi(t,x)$ any
smooth compactly supported test function, and we integrate on the
cylinder $Q=(0,T]\times B_{R}$. By definition of  local weak
solution, we obtain
\begin{equation}\label{Step.1.1}
\iint_Q \,\left[u^{p-1}\partial_t u
+\frac{4(p-1)}{(p+m-1)^2}\left|\nabla u^{\frac{p+m-1}{2}}\right|^2
\right]\,\psi^2\,\dx
\dt=-2\iint_Q\,u^{p+m-2}\nabla\psi\cdot\psi\nabla u \dx \dt.
\end{equation}
We now use Young's inequality: for any $\overrightarrow{a},
\,\overrightarrow{b}\in\RR^d$, and any $\delta>0$ we have
\[
|\overrightarrow{a} \cdot \overrightarrow{b}|
\le\frac{\delta}{2}|\overrightarrow{a}|^2+\frac{1}{2\delta}|\overrightarrow{b}|^2
\]
Together with H\"older's inequality, this  allows to estimate the
right-hand side of \eqref{Step.1.1}:
\begin{equation*}\label{Step.1.2}\begin{split}
-2\iint_Q\,u^{p+m-2} & \nabla\psi\cdot\psi\nabla u \dx \dt
=-\frac{4}{p+m-1}\iint_Q\,u^{\frac{p+m-1}{2}}\nabla\psi\cdot\psi\nabla u^{\frac{p+m-1}{2}} \dx \dt\\
&\le\frac{4}{p+m-1}\left[\frac{1}{2\delta}\iint_Q\,u^{p+m-1}\left|\nabla\psi\right|^2\dx
\dt
+\frac{\delta}{2}\iint_Q\,\left|\nabla u^{\frac{p+m-1}{2}}\right|^2\,\psi^2 \dx \dt\right]\\
&=\frac{2}{p-1}\iint_Q\,u^{p+m-1}\left|\nabla\psi\right|^2\dx \dt
+\frac{2(p-1)}{p+m-1}\iint_Q\,\left|\nabla u^{\frac{p+m-1}{2}}\right|^2\,\psi^2 \dx \dt\\
\end{split}
\end{equation*}
where in the last step we have chosen
$\delta=\frac{p-1}{p+m-1}>0$. Putting this calculation into
\eqref{Step.1.1}, we obtain
\begin{equation*}\label{Step.1.3}
\iint_Q \,u^{p-1}\partial_t u \,\psi^2\dx \dt
+\frac{2(p-1)}{(p+m-1)^2}\iint_Q \,\left|\nabla
u^{\frac{p+m-1}{2}}\right|^2 \,\psi^2\,\dx
\dt\le\frac{2}{p-1}\iint_Q\,u^{p+m-1}\left|\nabla\psi\right|^2\dx \dt
\end{equation*}
Now, we integrate the first term by parts (with respect to the
time variable)
\begin{equation*}\begin{split}
\iint_Q \,u^{p-1} & \partial_t u  \,\psi^2\dx \dt
=\frac{1}{p}\int_{B_R}\int_0^T
\,\partial_t(u^p)\,\psi^2\dx \dt\\
&=\frac{1}{p}\left[\int_{B_R}u^p(T,x)\psi^2(T,x)
\dx-\int_{B_R}u^p(0,x)\psi^2(0,x) \dx\right]-\frac{1}{p}\iint_Q
\,u^{p}\,\partial_t(\psi^2)\dx \dt\\
&=\frac{1}{p}\left[\int_{B_R}u^p(T,x)\psi^2(T,x)
\dx-\int_{B_R}u^p(0,x)\psi^2(0,x) \dx\right]-\frac{2}{p}\iint_Q
\,u^{p}\,\psi\,\partial_t(\psi)\,\dx \dt.
\end{split}
\end{equation*}
Collecting all the previous calculations, we obtain the first
basic inequality:
\begin{equation}\label{Step.1.FBI}\begin{split}
&\frac{1}{p}\left[\int_{B_R}u^p(T,x)\psi^2(T,x)
\dx-\int_{B_R}u^p(0,x)\psi^2(0,x) \dx\right]-\frac{2}{p}\iint_Q
\,u^{p}\,\psi\,\partial_t(\psi)\,\dx \dt\\
&+\frac{2(p-1)}{(p+m-1)^2}\iint_Q \,\left|\nabla
u^{\frac{p+m-1}{2}}\right|^2 \,\psi^2\,\dx
\dt\le\frac{2}{p-1}\iint_Q\,u^{p+m-1}\left|\nabla\psi\right|^2\dx \dt
\end{split}
\end{equation}
\noindent (ii) In order to continue, we assume that the test
function $\psi$ satisfies
\begin{itemize}
\item $0\le\psi(t,x)\le 1$, for any $(t,x)\in Q$, and
$\psi(0,x)=0$, for any $x\in B_R$
\item $\psi\equiv 1$ on $Q_1=[T_1,T]\times B_{R_1}\subset Q$ and
$\psi\equiv 0$ outside $Q$. Of course, we take $0\le R_1<R$ and
$0\le T_0<T_1\le T$.
\item Moreover, on $R\setminus R_1$, we assume that
\[
|\nabla\psi|\le\frac{c_\psi}{R-R_1}\qquad\mbox{and }\qquad
|\partial_t\psi|\le\frac{c_\psi^2}{T_1-T_0}\;.
\]
\end{itemize}
We may then write \eqref{Step.1.FBI} in the form
\[\begin{split}
\frac{p-1}{p}\int_{B_R}u^p(0,x)\psi^2(T,x) \dx
&+\frac{2(p-1)^2}{(p+m-1)^2}\iint_Q \,\left|\nabla u^{\frac{p+m-1}{2}}\right|^2\,\psi^2\,\dx \dt\\
&\le\,2\left[\iint_Q\,u^{p+m-1}\left|\nabla\psi\right|^2\dx \dt
+\frac{p-1}{p}\iint_Q
\,u^{p}\,\psi\,\left|\partial_t(\psi)\right|\,\dx \dt\right]\\
&\le 2c_\psi^2\,\left[
\frac{1}{(R-R_1)^2}+\frac{1}{T_1-T_0}\right]\left[\iint_Q\,\left(u^{p+m-1}+\,u^{p}\,\right)\dx
\dt
\right]\\
\end{split}
\]
We observe that
\[
\iint_{Q_1} \,\left|\nabla u^{\frac{p+m-1}{2}}\right|^2\,\dx
\dt\le\iint_Q \,\left|\nabla
u^{\frac{p+m-1}{2}}\right|^2\,\psi^2\,\dx \dt
\]
since $Q_1\subset Q$ and $\psi\equiv 1$ on $Q_1$, so that we
finally obtain
\[
\begin{split}
\mathcal{C}_{m,p}&\left[\int_{B_{R_1}}u^p(0,x) \dx +\iint_{Q_1}
\,\left|\nabla u^{\frac{p+m-1}{2}}\right|^2\,\dx \dt\right] \\&\le
2c_\psi^2\,\left[ \frac{1}{(R-R_1)^2}+\frac{1}{T_1-T_0}\right]
\left[\iint_Q\,\left(u^{p+m-1}+\,u^{p}\,\right)\dx \dt\right],
\end{split}
\]
where
\begin{equation}\label{form.cmp}  \mathcal{C}_{m,p}=\min\left\{
\frac{p-1}{p}\,,\frac{2(p-1)^2}{(p+m-1)^2}\right\}.
\end{equation}
As a conclusion, we have obtained \eqref{Step.1.final} with
precise constant
\begin{equation}\mathcal{C}(m,p) =
2c_\psi^2\,\mathcal{C}_{m,p}^{-1}\;.
\end{equation}
Note that $\mathcal{C}_{m,p}$ depends
also on $d$ though we are not indicating it.

 We conclude  by noticing that the proof can be repeated
for $u$ sub-solution (with the same regularity), that means, when
$u_t\le \Delta u^m$  the above estimate continues to hold.\qed

\medskip

\noindent {\bf Improving the constant.} We would like to eliminate
the dependence of $\mathcal{C}(m,p)$ on $p$ in what follows since
$p$ will vary (in an increasing way). This dependence takes place
through $\mathcal{C}_{m,p}$. Now, for $m\ge 0$ it is easy to see
that $\mathcal{C}_{m,p}\ge (p-1)/p$ \ and we have to assume that
$p\ge p_0>1$, so that $\mathcal{C}(m,p)$ is bounded by an
expression that depends only on $p_0 $ and $d$.

For $m<0$, since we have $p>1-m$ we get $(p-1)/p> |m|/(1-m)$. A
lower bound for $\mathcal{C}_{m,p}$ needs the last term to be
bounded above, and this implies that $p$ must be away from $1-m$,
so that we assume that $p\ge p_0'=(1+\alpha) (1-m)$ for some
$\alpha>0$ in which case we get
\begin{equation}\label{Step.1.C0m}
\mathcal{C}_{m,p}\ge \min\left\{\dfrac{|m|}{1-m}, \
2\left(1+\dfrac{|m|}{\alpha (1-m)}\right)^2\right\}\,:=\,\mathcal{C}(m).
\end{equation}
In any case we may write $\mathcal{C}(m)$ instead of
$\mathcal{C}_{m,p}$ if the family of $p$'s fulfills the stated
conditions.

\medskip

\noindent {\bf Final result of Step 1.} We need to improve Lemma
\ref{Lem.Step.1} in the following way
\begin{cor}
Under the running assumptions, for every $m<1$ and
$p>\max\{1,1-m\},$ then for any $T_0<T_1<T$, $0<R_1<R$ we have
\begin{equation}\label{Step.1.final.sup}
\begin{split}
\sup_{s\in(T_1,T)}\int_{B_{R_1}}u^p(s,x) \dx
&+\int_{T_1}^{T}\int_{B_{R_1}} \,\left|\nabla
u^{\frac{p+m-1}{2}}\right|^2\,\dx \dt\\
&\le \mathcal{C}(m)\,\left[
\frac{1}{(R-R_1)^2}+\frac{1}{T_1-T_0}\right]
\left[\int_{T_0}^{T}\int_{B_R}\,\left(u^{p+m-1}+\,u^{p}\,\right)\dx
\dt\right].
\end{split}
\end{equation}
Moreover, if $u$ is a sub-solution, and $u\ge 1$, we have
\begin{equation}\label{Step.1.u>1}
\begin{split}
\sup_{s\in(T_1,T)}\int_{B_{R_1}}u^p(s,x) \dx
&+\int_{T_1}^{T}\int_{B_{R_1}} \,\left|\nabla
u^{\frac{p+m-1}{2}}\right|^2\,\dx \dt
\le \mathcal{C}(m)\,\left[
\frac{1}{(R-R_1)^2}+\frac{1}{T_1-T_0}\right]
\left[\int_{T_0}^{T}\int_{B_R}\,u^{p}\,\dx \dt\right]\\
\end{split}
\end{equation}
with $\mathcal{C}(m)$ as in \eqref{Step.1.C0m}.
\end{cor}

\noindent {\sl Proof.~} First we recall a property of the supremum:
there exists a $t_0\in
(T_1,T]$ such that
\[
\frac{1}{2}\sup_{s\in(T_1,T)}\int_{B_{R_1}}u^p(s,x) \dx
    \le \int_{B_{R_1}}u^p(t_0,x) \dx
\]
We use this observation together to the result of Lemma \ref{Lem.Step.1}, in two different ways:

 \noindent (i) We use Lemma \ref{Lem.Step.1} with the
$T$ replaced by $t_0$ and still keeping   $0\le T_0<T_1<t_0$. We
get
\begin{equation}\label{Step.1.Sup.1}
\begin{split}
\frac{1}{2}\sup_{s\in(T_1,T)}\int_{B_{R_1}}u^p(s,x) \dx
    &\le\int_{B_{R_1}}u^p(t_0,x) \dx\\
    &\le \mathcal{C}(m)\,\left[ \frac{1}{(R-R_1)^2}+\frac{1}{T_1-T_0}\right]
        \left[\int_{T_0}^{t_0}\int_{B_R}\,\left(u^{p+m-1}+\,u^{p}\,\right)\dx \dt\right]\\
\end{split}
\end{equation}
In the sequel recall that $t_0\le T$.

 \noindent (ii) Next, we choose the same $T_1$ and we apply Lemma \ref{Lem.Step.1}, to get
\begin{equation}\label{Step.1.Sup.2}
\begin{split}
\int_{T_1}^{T}\int_{B_{R_1}}\,\left|\nabla u^{\frac{p+m-1}{2}}\right|^2\,\dx \dt
    \le \mathcal{C}(m)\,\left[ \frac{1}{(R-R_1)^2}+\frac{1}{T_1-T_0}\right]
        \left[\int_{T_0}^{T}\int_{B_R}\,\left(u^{p+m-1}+\,u^{p}\,\right)\dx \dt\right].
\end{split}
\end{equation}

\noindent Summing up the two inequalities \eqref{Step.1.Sup.1} and
\eqref{Step.1.Sup.2} gives  the desired inequality
\eqref{Step.1.final.sup}\,.

\noindent For the last part, we remark that if we apply inequality
\eqref{Step.1.final.sup} to a sub-solution $u\ge 1$, then
$u^{p+m-1}\le u^p$ so that we obtain \eqref{Step.1.u>1}\,, and the
proof is thus concluded.\qed

\subsubsection*{ Step 2. Iterative form of the Sobolev
Inequality}

The next lemma is just a different form of the usual Sobolev
inequality, adapted to our aims.

\begin{lem}\label{Lem.Sob.Iter} Let $f\in \LL^2(Q)$ with $\nabla f\in \LL^2(Q)$. We then have
\begin{equation}\label{Step.2.Sob.Iter}\begin{split}
\int_{T_1}^T\,\int_{B_R}f^{2\sigma}\dx \dt
    &\le\,2\mathcal{S}_2^2
        \left[\int_{T_1}^T\,\int_{B_R}\left(f^2 + R^2\,\big|\nabla f\big|^2\right)\dx\dt\right]
    \,\sup_{s\in(T_1,T)}\left[\frac{1}{R^d}\int_{B_R}f^{2(\sigma-1)q}(s,x)\dx\right]^{\frac{1}{q}}
\end{split}
\end{equation}
for any $\sigma\in\big(1,\sigma^*\big)$, and for any $0\le T_1< T$ and $R>0$\,, where
\begin{equation}
\sigma^* =\frac{2^*}{2}=\left\{\begin{array}{lll}
\frac{d}{d-2}=\frac{1}{m_c}\;,&\mbox{if~}d\ge 3\\
2\;,&\mbox{if~}d=1,2\\
\end{array}\right.\qquad\mbox{and}\qquad
q =\frac{\sigma^*}{\sigma^*-1}=\left\{\begin{array}{lll}
\frac{d}{2}\;,&\mbox{if~}d\ge 3\\
2\;,&\mbox{if~}d=1,2.\\
\end{array}\right.
\end{equation}
Here,  $\mathcal{S}_2$ is the constant of the classical Sobolev
inequality $ \|f\|_{2^*}\le\mathcal{S}_2\,\left(\|\nabla f\|_2 +
\|f\|_2\right) $,  with $2^*=2d/(d-2)$ for $d\ge 3$ and $2^*=4$
for $d=1,2$.
\end{lem}
\noindent {\sl Proof.~} Since the estimate \eqref{Step.2.Sob.Iter}
is scaling invariant, it is sufficient to prove it for $R=T-T_1=1$,
and we denote by $B=B_1$ the unit ball of $\RR^d$. By Sobolev
 and H\"older inequalities we then get
\[
\begin{split}
\int_B f^{2\sigma}\dx&=\int_B f^2
f^{2(\sigma-1)}\dx\le\left[\int_B f^{2^*}
\dx\right]^{\frac{2}{2^*}}\left[\int_B
f^{2(\sigma-1)q}\dx\right]^{\frac{1}{q}}\\
&\le2\mathcal{S}_2^2\left[\int_B\big|\nabla f\big|^2\dx+\int_B
f^2\dx
\right]\sup_{s\in(0,1)}\left[\int_{B_1}f^{2(\sigma-1)q}(s,x)\dx\right]^{\frac{1}{q}}
\end{split}
\]
Integrating in time over $(0,1)$ and rescaling back, gives
inequality \eqref{Step.2.Sob.Iter}.\qed

\subsubsection*{Step 3. The Iteration } In this step we use
the inequalities of the preceding steps to start the iteration in
the Moser style. We first define
$v(t,x)=\max\big\{u(t,x),1\big\}$. Then we observe that when $u$
is a local weak solution to $u_t=\Delta\,u^m$\,, then $v$ is
a local weak {\sl sub-solution} to $v_t=\Delta\,v^m$. It is clear
that $u\le v\le 1+u$ for almost any $(t,x)\in Q$\,.

\medskip

\noindent \textsc{Preparation of the iteration step.} Letting
$f^2=v^{p+m-1}$ in the modified Sobolev inequality
\eqref{Step.2.Sob.Iter} gives
\begin{equation}\label{Step.3.1}
\begin{split}
\iint_{Q_1} v^{\sigma(p+m-1)}\dx \dt \le 2 &\mathcal{S}_2^2
\left[\iint_{Q_1}\left(v^{p+m-1} +
R_1^2\,\big|\nabla v^{\frac{p+m-1}{2}}\big|^2\right)\dx \dt\right]
\left[\sup_{t\in(T_1,T)}\frac{1}{R_1^d}\int_{B_{R_1}}v^{(p+m-1)(\sigma-1)q}\dx
\right]^{\frac{1}{q}}
\end{split}
\end{equation}
where $Q_1=(T_1, T]\times B_{R_1}\subset Q_0=(T_0,T]\times B_{R}$.

\noindent Since we have assumed that $v\ge 1$, then $v^{p+m-1}\le
v^p$ and we can use the energy inequality  \eqref{Step.1.u>1} to
estimate the two terms of the right hand side of the above
inequality \eqref{Step.3.1}, in terms of the same quantity. First
we estimate
\[
\begin{split}
\iint_{Q_1}\left(v^{p+m-1} + R_1^2\,\big|\nabla
v^{\frac{p+m-1}{2}}\big|^2\right) \dx \dt
    &\le \iint_{Q_1}v^{p}\dx \dt\\
    &+R_1^2\,\mathcal{C}(m)\,\left[
        \frac{1}{(R_0-R_1)^2}+\frac{1}{T_1-T_0}\right]
        \left[\iint_{Q_0}\,v^{p}\,\dx \dt\right]\\
    &\le 2\,R_1^2\,\mathcal{C}(m)\,
        \left[\frac{1}{(R_0-R_1)^2}+\frac{1}{T_1-T_0}\right]
        \left[\iint_{Q_0}\,v^{p}\,\dx \dt\right].
\end{split}
\]
In the last step we use the fact that
$R_1^2\,\,\mathcal{C}(m)\,\left[
\frac{1}{(R_0-R_1)^2}+\frac{1}{T_1-T_0}\right]\ge 1\;, $ which is not
restrictive.

We next  estimate the sup term, again using the energy inequality
\eqref{Step.1.u>1}\,, but we replace $p$ with
$(p+m-1)(\sigma-1)q$. If the exponent is larger than
$\max\{1,1-m\}$  we obtain
\[
\sup_{t\in(T_1,T)}\frac{1}{R_1^d}\int_{B_{R_1}}v^{(p+m-1)(\sigma-1)q}\dx
    \le \frac{\mathcal{C}(m)}{R_1^d}\,\left[
        \frac{1}{(R-R_1)^2}+\frac{1}{T_1-T_0}\right]
        \left[\iint_{Q_0}\,v^{(p+m-1)(\sigma-1)q}\,\dx \dt\right]\\
\]
Summing up, we have estimated \eqref{Step.3.1} as follows
\begin{equation}\label{Step.3.2}
\begin{split}
\iint_{Q_1} v^{\sigma(p+m-1)}\dx \dt
    &\le 4 \mathcal{S}_2^2\;\,R_1^{2-\frac{d}{q}}\,\mathcal{C}(m)^{1+\frac{1}{q}}\,
        \left[\frac{1}{(R_0-R_1)^2}+\frac{1}{T_1-T_0}\right]^{1+\frac{1}{q}}\\
    &\times\left[\iint_{Q_0}\,v^{p}\,\dx \dt\right]
        \left[\iint_{Q_0}\,v^{(p+m-1)(\sigma-1)q}\,\dx \dt\right]^{\frac{1}{q}}
\end{split}
\end{equation}
Finally, we remark that $R_1^{2-\frac{d}{q}}=1$, since
$\frac{d}{q}=2$, $q$ being defined as in Lemma
\ref{Lem.Sob.Iter}\,.

\medskip

\noindent \textsc{The First Iteration Step.} We now use
\eqref{Step.3.2} in the following way: first we choose
$\sigma\in(1,\sigma^*)$, where $\sigma^*$ is as in Lemma
\ref{Lem.Sob.Iter}, in such a way that
\[
(p+m-1)(\sigma-1)q=p\qquad\mbox{that is}\qquad
\sigma=1+\frac{p}{q(p+m-1)}.
\]
A straightforward calculation shows that $\sigma\in(1,\sigma^*)$
if and only if $p>p_c$\,. This is where the restriction on $p$
appears for the first time.

We are now ready to begin with the first iterative step, by
letting
\[
p_0=p=(p+m-1)(\sigma-1)q\,,\qquad\mbox{and }\qquad
p_1=(p_0+m-1)\sigma=p_0\left(1+\frac{1}{q}\right)+m-1\,.
\]
We remark that
\[
p_1>p_0 \iff p_0>p_c=\frac{d(1-m)}{2}\,.
\]
Estimate \eqref{Step.3.2} now becomes
\begin{equation}\label{Step.3.Iter.01}
\begin{split}
\iint_{Q_1}v^{p_1}\dx \dt
    &\le 4 \mathcal{S}_2^2\;\,\mathcal{C}(m)^{1+\frac{1}{q}}\,
        \left[\frac{1}{(R_0-R_1)^2}+\frac{1}{T_1-T_0}\right]^{1+\frac{1}{q}}
        \left[\iint_{Q_0}\,v^{p_0}\,\dx \dt\right]^{1+\frac{1}{q}}\\
    &= I_{0,1}\left[\iint_{Q_0}\,v^{p_0}\,\dx \dt\right]^{1+\frac{1}{q}}
\end{split}
\end{equation}
which is the first iterative step.

\medskip

\noindent\textsc{The $k$-th Iteration Step.} Letting
\[
p_{k+1}=p_k\left(1+\frac{1}{q}\right)+m-1\;,\qquad\mbox{with}\qquad
p_{k+1}>p_k \iff p_k\ge p_0>p_c,
\]
we get the iterative inequality
\begin{equation}\begin{split}
\left[\iint_{Q_{k+1}}v^{p_{k+1}}\dx
\dt\right]^{\frac{1}{p_{k+1}}}&\le
I_{k,k+1}^{\frac{1}{p_{k+1}}}\left[\iint_{Q_k}\,v^{p_k}\,\dx
\dt\right]^{\frac{1}{p_k}\left(1+\frac{1}{q}\right){\frac{p_k}{p_{k+1}}}}.
\end{split}
\end{equation}
In order to find a convenient value for $I_{k,k+1}$ we choose a
decreasing sequence of radii $R_\infty \longleftarrow
R_{k+1}<R_k<R_0$ such that $0<R_k-R_{k+1}=\rho/k^2$, and  a
sequence of times $0\le T_0\le T_k \le T_{k+1}\to T_\infty<T$ such
that $T_{k+1}-T_k=\tau/k^4$. This means taking
\begin{equation}\label{rhotau}
\rho=c_1(R_0-R_\infty), \qquad \tau=c_2(T_\infty-T_0),
\end{equation}
 with $c_1=1/\left(\sum_{k=0}^{+\infty}k^{-2}\right)>0$, and
$c_2=1/\left(\sum_{k=0}^{+\infty}k^{-4}\right)>0$. Then,
\begin{equation}\label{Iter.est.Ik}
\begin{split}
I_{k,k+1}
    &=4 \mathcal{S}_2^2\,\mathcal{C}(m)^{1+\frac{1}{q}\,}\,
        \left[\frac{1}{(R_{k}-R_{k+1})^2}+\frac{1}{T_{k+1}-T_k}\right]^{1+\frac{1}{q}}
  \le 4\mathcal{S}_2^2\,\mathcal{C}(m)^{1+\frac{1}{q}}\,
        \,\left(2\left(\rho^{-2}+\tau^{-1}\right)\,k^4\right)^{1+\frac{1}{q}}\\
    &\le\,\left[2\mathcal{S}_2\right]^2\,
    \left[2\left(\rho^{-2}+\tau^{-1}\right)\, \mathcal{C}(m)\right]^{1+\frac{1}{q}}\,
        \,\left(k^4\right)^{1+\frac{1}{q}}
     = J_0\,J_{1}^{1+\frac{1}{q}}\,\left(k^4\right)^{1+\frac{1}{q}}\\
\end{split}
\end{equation}
 We now calculate the exponents $p_k$:
\begin{equation}\label{pk}
\begin{split}
p_{k+1}&=p_k\left(1+\frac{1}{q}\right)+m-1
=\left[1+\frac{1}{q}\right]^{k+1}p_0+(m-1)\sum_{n=0}^{k}\left[1+\frac{1}{q}\right]^{n}\\
&=\left[1+\frac{1}{q}\right]^{k+1}\left[p_0+(m-1)\sum_{j=1}^{k+1}\left[1+\frac{1}{q}\right]^{-j}\right]
=\big[p_0-q(1-m)\big]\left[1+\frac{1}{q}\right]^{k+1}+q(1-m)\,.\\
\end{split}
\end{equation}
notice that
\[
\lim_{k\to\infty}\frac{\left[1+\frac{1}{q}\right]^{k+1}}{p_{k+1}}=\frac{1}{p_0+q(m-1)}
\qquad\mbox{and}\qquad
\lim_{k\to\infty}\frac{1}{p_{k+1}}\sum_{j=0}^k\left[1+\frac{1}{q}\right]^{j}=\frac{q}{p_0+q(m-1)}\,.
\]
The iterative step now reads
\begin{equation}\label{Iter.step.k+1}
\begin{split}
\left[\iint_{Q_{k+1}}v^{p_{k+1}}\dx \dt\right]^{\frac{1}{p_{k+1}}}
    &\le\,I_{k,k+1}^{\frac{1}{p_{k+1}}}\;\;I_{k-1,k}^{\left[1+\frac{1}{q}\right]\frac{1}{p_{k+1}}}\;\;
        \ldots\;\;I_{0,1}^{\left[1+\frac{1}{q}\right]^{k}\frac{1}{p_{k+1}}}
        \,\left[\iint_{Q_{0}}v^{p_{0}}\dx\dt\right]^{\frac{\left[1+\frac{1}{q}\right]^{k+1}}{p_{k+1}}}\\
\end{split}
\end{equation}
Now we use \eqref{Iter.est.Ik} to estimate
\begin{equation}
\begin{split}
I_{k,k+1}^{\frac{1}{p_{k+1}}}\;\;I_{k-1,k}^{\left[1+\frac{1}{q}\right]
        \frac{1}{p_{k+1}}}\;\;\ldots\;\;I_{0,1}^{\left[1+\frac{1}{q}\right]^{k}\frac{1}{p_{k+1}}}
    &\le\left[J_0\,J_1^{1+\frac{1}{q}}\right]^{\frac{1}{p_{k+1}}\sum_{j=0}^k\left[1+\frac{1}{q}\right]^{j}}\\
    &\times\;k^{4\frac{1}{p_{k+1}}}\,(k-1)^{4\frac{1+\frac{1}{q}}{p_{k+1}}}
        \,(k-2)^{4\frac{(1+\frac{1}{q})^2}{p_{k+1}}}
        \;\;\ldots\;\;2^{4\frac{(1+\frac{1}{q})^{k-2}}{p_{k+1}}}\,1\\
    &=\left[J_0\,J_1^{1+\frac{1}{q}}\right]^{\frac{1}{p_{k+1}}\sum_{j=0}^k\left[1+\frac{1}{q}\right]^{j}}
        \,\prod_{j=1}^{k} j^{4\frac{(1+\frac{1}{q})^{k-j}}{p_{k+1}}}\\
\end{split}
\end{equation}
Moreover, passing to the limit in \eqref{Iter.step.k+1} when
$k\to\infty$\,, we get (we refer to Appendix A3 for further details)
\begin{equation}
\sup_{Q_\infty}|v|\,\le\,J_0^{\frac{q}{p_0+q(m-1)}}\,J_1^{\frac{q+1}{p_0+q(m-1)}}\,s_1\ee^{4(q+1)}
\;\left[\iint_{Q_{0}}v^{p_{0}}\dx
\\dt\right]^{\frac{1}{p_0+q(m-1)}}.
\end{equation}
We have estimated the constants to ensure that they remain bounded
in the limit $k\to +\infty$, see the Appendix for the details.
Moreover, these constants blow up when $R_\infty\to R_0$ or
$T_\infty \to T_0$: indeed, while $J_0$ and $s_1$ only depend on
$m,d$ and $p$, the constant $J_1$ depends on $\rho$ and $\tau$:
and blows up as $T_0-T_\infty\to 0$ or  $R_0-R_\infty\to 0$ since
$J_1\sim\left(\rho^{-2}+\tau^{-1}\right)$. We have finished this
part of the iteration since $q=d/2$ which gives in the sequel the
correct exponent in the last integral. In case $d=1,2$ we have to
observe that $q=2$ so that the exponent is $s=1/(p_0+2(m-1)$.

\noindent This fact forces the final cylinder $Q_\infty=(T_\infty,
T]\times B_{R_\infty}\subset Q_0=(T_0, T]\times B_{R_0}$ to be
strictly contained in the initial one. We  obtain
\begin{equation}\label{sup.est.mean.time}
\sup_{Q_\infty}|v|\,\le\,\mathcal{C}_{\rm
loc}\,\left[\frac{1}{(R_0-R_\infty)^2}+\frac{1}{T_\infty-T_0}\right]^{\frac{q+1}{p_0+q(m-1)}}
\;\left[\iint_{Q_{0}}v^{p_{0}}\dx
\dt\right]^{\frac{1}{p_0+q(m-1)}}\,.
\end{equation}
\noindent Notice that $\mathcal{C}_{\rm loc}$ only depends on
$m,d$ and $p$. We conclude the proof by going back from $v$ to
$u$, using the fact that $u\le v\le u+1$, by definition of $v$.
>From \eqref{sup.est.mean.time} we easily get
\[
\begin{split}
\sup_{Q_\infty}|u|\,
    \le\sup_{Q_\infty}|v|\,
    &\le\,\mathcal{C}_{\rm loc}\,
        \left[\frac{1}{(R_0-R_\infty)^2}+\frac{1}{T_\infty-T_0}\right]^{\frac{q+1}{p_0+q(m-1)}}
        \;\left[\iint_{Q_{0}}v^{p_{0}}\dx \dt\right]^{\frac{1}{p_0+q(m-1)}}\\
    &\le\,\mathcal{C}_{\rm loc}\,
        \left[\frac{1}{(R_0-R_\infty)^2}+\frac{1}{T_\infty-T_0}\right]^{\frac{q+1}{p_0+q(m-1)}}
        \;\left[\iint_{Q_{0}}u^{p_{0}}\dx \dt+\vol\big(Q_0\big)\right]^{\frac{1}{p_0+q(m-1)}}\,.
\end{split}
\]
This concludes the proof, after changing the notation, putting
$p=p_0$, $R_\infty=R_1<R_0$\,, $T_\infty=T_1>T_0$, and $q=d/2$
since we are dealing with $d\ge 3$\,.\qed

\subsection{\large Local Smoothing effect. Proof of
Theorem~\ref{thm.upper}} \label{sect4}

\noindent We now combine the results of
Theorem~\ref{Lem.Local.Lp.norms}  and of
Theorem~\ref{Thm.1.Integrated}, to prove the local smoothing
effect in the form described in Theorem~\ref{thm.upper}. We
consider $u$ defined in $(0,T)\times B_{R_0}(x_0)$, then
take a smaller radius  $R_1$  and write
$R_1=(1-\varepsilon)\varrho $ and $R_0=(1+\varepsilon)\varrho$:
this defines $\rho$ and $\varepsilon$. We then consider the rescaled
solution
\begin{equation}
\widehat{u}(t,x)=K\,u\left(\tau\,t\,,\varrho\,x+x_0\right)\,,\qquad
K=\left(\frac{\varrho^2}{\tau}\right)^{\frac{1}{1-m}}
\end{equation}
with $0<\tau<T$. Then, we apply to $\widehat{u}$ the result of
Theorem \ref{Thm.1.Integrated} over the cylinders $Q_0=(0,1]\times
B_1$ and $Q_1=(\varepsilon^2,1]\times B_{1-\varepsilon}$, for some
$\varepsilon\in(0,1)$\,, so that the bound reads
\[\begin{split}
\sup_{Q_1}|\widehat{u}|\,
    &\le\,\frac{\mathcal{C}_{\rm loc}}{\varepsilon^{2\frac{q+1}{p+q(m-1)}}}
        \;\left[\iint_{Q_{0}}\widehat{u}^{p}\dx \dt+\vol\big(Q_0\big)\right]^{\frac{1}{p+q(m-1)}}
    \le\,\frac{\mathcal{C}_{\rm loc}}{\varepsilon^{2\frac{q+1}{p+q(m-1)}}}
        \;\left[\iint_{Q_{0}}\widehat{u}^{p}\dx \dt+\omega_d\right]^{\frac{1}{p+q(m-1)}}
\end{split}
\]
since $\vol\big(Q_0\big)=\omega_d$\,. Moreover, we know that
$\mathcal{C}_{\rm loc}$ only depends on $m,d,p$\,.

\noindent Next, we estimate the time integral, using Theorem
\ref{Lem.Local.Lp.norms} applied to the rescaled solution
$\widehat{u}$ on the balls $B_{1}\subset B_{1+\varepsilon}$ and
for times $t\in[0,1]$ and for $p>p_c>1-m$\,:
\begin{equation}
\int_{B_1}|\widehat{u}(t,x)|^{p} \dx\le
2^{\frac{p}{1-m}-1}\int_{B_{1+\varepsilon}}|\widehat{u}(0,x)|^{p}\dx +
2^{\frac{p}{1-m}-1}\left[K_{\varepsilon,
p}\;t\right]^{\frac{p}{1-m}},
\end{equation}
where
\[
K_{\varepsilon, p}=
\frac{p\,c_{m,d}}{\varepsilon^{2}}\,\vol\left(B_{1+\varepsilon
}\setminus B_1\right)^{(1-m)/p}\le p\,c'_{m,d}\varepsilon^{\frac{1-m}{p}-2}
\,
\]
An integration in time over $(0,1)$ gives
\begin{equation*}\begin{split}
\iint_{Q_{0}}\widehat{u}^{p_{0}}\dx
\dt=\int_0^1\int_{B_1}\widehat{u}^{p_{0}}\dx \dt&\le
2^{\frac{{p}}{1-m}-1}\int_{B_{1+\varepsilon}}|\widehat{u}(0,x)|^{p}\dx
+\frac{2^{\frac{{p}}{1-m}-1}}{\frac{{p}}{1-m}+1}\left[p\,c'_{m,d}\varepsilon^{\frac{1-m}{p}-2}\right]^{\frac{p}{1-m}}\\
&=\mathcal{S}_0\int_{B_{\lambda}}|\widehat{u}(0,x)|^{p}\dx
+\frac{\mathcal{K}_{m,d,p}}{\varepsilon^{\frac{2p}{(1-m)}-1}}
\end{split}
\end{equation*}
and we remark that $\mathcal{K}_{m,d,p}$ and $\mathcal{S}_0$ only depend on $m,p$\,.

\noindent Now we put together the above estimates, and we rescale
back from $\widehat{u}$ to $u$, also changing variable $\varrho
x=y$ in the integrals, and we obtain (using
$K^{1-m}=\varrho^2/\tau$)
\begin{equation}\begin{split}
\sup_{(s,y)\in(\varepsilon^2\tau,\tau]\times
B_{(1-\varepsilon)\varrho}}u(s,y)
    &\le\,K^{-1}\frac{\mathcal{C}_{\rm loc}}{\varepsilon^{2\frac{q+1}{p+q(m-1)}}}
        \;\left[\mathcal{S}_0\int_{B_{1+\varepsilon}}|\widehat{u}(0,x)|^{p}\dx
        +\frac{\mathcal{K}_{m,d,p}}{\varepsilon^{\frac{2p}{(1-m)}-1}}+
        \omega_d\right]^{\frac{1}{p+q(m-1)}}\\
    &\le\,\varrho^{\frac{2q-d}{p+q(m-1)}} \frac{\mathcal{C}_{\rm loc}}
    {\varepsilon^{2\frac{q+1}{p+q(m-1)}}}
        \frac{\kappa_1\,\mathcal{S}_0^{\frac{1}{p+q(m-1)}}}{\tau^\frac{q}{p+q(m-1)}}
        \left[\int_{B_{(1+\varepsilon)\varrho}}|u(0,y)|^{p}\rd y \right]^{\frac{1}{p+q(m-1)}}\\
    &+\kappa_2\,\,\frac{\mathcal{C}_{\rm loc}}{\varepsilon^{2\frac{q+1}{p+q(m-1)}}}
        \left[\,\frac{\mathcal{K}_{m,d,p}}{\varepsilon^{\frac{2p}{(1-m)}-1}}+\omega_d\Big)\right]^{\frac{1}{p+q(m-1)}}
        \left[\frac{\tau}{\varrho^2}\right]^{\frac{1}{1-m}}.
\end{split}
\end{equation}
In dimension $d\ge 3$ we have taken $q=d/2$ which allows to cancel
the appearance of $\varrho$ in the first term of the right-hand
side and simplify the dependence on $\tau$. We then have
\begin{equation}\label{bit.more.sup}
\sup_{(s,y)\in(\varepsilon^2\tau,\tau]\times B_{R_1}}u(s,y)\le
\frac{\overline{\mathcal{C}}_{1}}{\tau^{d\vartheta_p}}
        \,\left[\int_{B_{R_0}}|u_0(x)|^{p}\dx\right]^{2\vartheta_p}
        +\overline{\mathcal{C}}_{2}\left[\frac{\tau}{\varrho^2}\right]^{\frac{1}{1-m}}
\end{equation}
where we have also  used $(a+b)^\sigma\le\kappa_1 a^\sigma
+\kappa_2 b^{\sigma}$. Putting $\tau=t$ we have obtained in
particular the desired formula \eqref{upper.est}\,. We last remark
that $\varepsilon$ must be strictly positive even if it can be
chosen arbitrarily small; in the limit $\varepsilon\to 0$, the
constants $\overline{\mathcal{C}}_i$ blow up, since
\begin{equation}\label{Const.upper}
\overline{\mathcal{C}}_1=\kappa_1\,\mathcal{S}_0^{\frac{1}{p+q(m-1)}}
        \frac{\overline{\mathcal{C}}_{\rm loc}}{\varepsilon^{2\frac{q+1}{p+q(m-1)}}}
    \qquad\mbox{and}\qquad
\overline{\mathcal{C}}_2=\kappa_2\,\,\frac{\overline{\mathcal{C}}_{\rm
loc}}{\varepsilon^{2\frac{q+1}{p+q(m-1)}}}
\left[\,\frac{\mathcal{K}_{m,d,p}}{\varepsilon^{\frac{2p}{(1-m)}-1}}+\omega_d\right]^{\frac{1}{p+q(m-1)}}
\end{equation}
We conclude the proof by switching to the same notations as in the
statement of Theorem~\ref{thm.upper}, just by substituting
$R_1=(1-\varepsilon)\varrho $ and $R_0=(1+\varepsilon)\varrho$\,,
it is clear that the result holds for any $R_1<R_0$, and that the
constants $\overline{\mathcal{C}}_i$ blow up when $R_1\to R_0$. We recall that $\mathcal{C}_{\rm loc}$ only depends on $m,d,p$\,.

We have thus proved the following result:

\begin{thm} 
Let $p\ge 1$ if $m>m_c$ or  $p>p_c$ if $m\le m_c$.
Then there are positive constants $\overline{\mathcal{C}}_{1}$,
$\overline{\mathcal{C}}_{2}$ such that for any $0<R_1<R_0$ we have
\begin{equation}\label{upper.est.old}
\begin{split}
\sup_{(s,y)\in(t_0,t]\times B_{R_1}}u(s,y)&\le\frac{\overline{\mathcal{C}}_{1}}{t^{d\vartheta_p}}\,\left[
\int_{B_{R_0}}|u_0(x)|^{p}\dx
\right]^{2\vartheta_p}+\overline{\mathcal{C}}_{2}\left[\frac{t}{R_0^2}\right]^{\frac{1}{1-m}}.
\end{split}
\end{equation}
where $t_0=\big[(R_0-R_1)/(2R_0)\big]^2\,t$ and the constants $\overline{\mathcal{C}}_{i}$ depend on $m,d$
and $p$, $R_1$ and $R_0$ and blow up when $R_1/R_0\to 1$; an
explicit formula for $\overline{\mathcal{C}}_{i}$ is given by \eqref{Const.upper}.
\end{thm}
The above theorem is nothing but a slightly stronger form for Theorem \ref{thm.upper}:
we just take the limit $R_1\to 0$ in inequality \eqref{upper.est.old} to obtain \eqref{upper.est}. The final expression for the constants $\mathcal{C}_i$ in \eqref{upper.est} corresponds to the limit of $\overline{\mathcal{C}}_i$ in \eqref{Const.upper} as $\varepsilon\to 1$ (i.e. $R_1\to 0$) and do not depend on the radii. This concludes the proof of Theorem \ref{thm.upper}.  \qed


\section{ \large  Part III. Harnack Inequalities}
\label{sect.parabolicHI}

By  joining together the local upper and lower estimates obtained in
Parts I and II we can draw interesting conclusions in terms of
special forms of Harnack Inequalities. These are expressions
relating the maximum and minimum of a solution inside certain
cylinders. In the standard case one has
 \begin{equation}
\sup_{Q_1} u(t,x)\le C \,\inf_{Q_2} u(t,x),
\end{equation}
see \cite{Moser} and \cite{TR}.  The main
idea is that  the formula applies for a large class of solutions
and the constant $C$ that enters the relation does not depend on
the particular solution, but only on the data like $m,d$ and the
size of the cylinder $R$, but not on time.
The cylinders in the standard case are supposed to be ordered
in time, $Q_1=[t_1,t_2]\times B_R(x_0)$,
$Q_2=[t_3,t_4]\times B_R(x_0)$, with $t_1\le t_2<t_3\le t_4$.

It is well-known that in the degenerate nonlinear elliptic or
parabolic problems a plain form of the inequality does not hold.
In the  work of DiBenedetto and collaborators, see the book
\cite{DBbook} or the recent work \cite{DGV}, versions are obtained
where some information of the solution is used to define so-called
intrinsic sizes, like the size of the parabolic cylinder(s), that
usually depends on $u(t_0,x_0)$. They are called {\sl intrinsic
Harnack inequalities}. The authors of \cite{DGV}  show that the
size of a convenient cylinder for the Harnack inequality to hold
has the form
\[
I_R(t_0,x_0)= \left(t_0-c\,u(t_0,x_0)^{1-m}R^2,
t_0+c\,u(t_0,x_0)^{1-m}R^2\right)\times B_R(x_0)
\]
with a  fixed constant  $c>0$ which depends only on $m,d$\,, that
can be chosen ``a priori'', but only in the good range
$m_c<m<1$\,.  This cylinder is called intrinsic because it depends
on the value of the solution $u$ at a given point $(t_0,x_0)$\,.

The Harnack Inequalities of \cite{DKV,DGV}, in the supercritical range then read:
\textit{There exist positive constants $\overline{c}$ and
$\overline{\delta}$ depending only on $m,d$, such that for all
$(t_0,x_0)\in Q=(0,T)\times\Omega$ and all cylinders of the type
$I_{8R}\subset Q$, we have
\[
\overline{c}\,u(t_0,x_0)\le \inf_{x\in B_R(x_0)}u(t,x)
\]
for all times $t_0-\overline{\delta}\,u(t_0,x_0)^{1-m}\,R^2<t<t_0
+\overline{\delta}\,u(t_0,x_0)^{1-m}\,R^2$. The constants
$\overline{\delta}$ and $\overline{c}$ tend to zero as
$m\to 1$ or as $m\to m_c$\,.}

They also give a counter-example in the lower range $m<m_c$, by
producing an explicit local solution that does not satisfy any kind
of Harnack inequality (neither of the types called intrinsic,
elliptic, forward,  backward) if one fixes ``a priori'' the
constant $c$\,. At this point a natural question is posed:

\textit{What form, if any, the Harnack estimate might take for $m$
in the sub–critical range $0<m\le m_c$?}

\noindent The following is an answer this question.

\medskip

\noindent {\sc New approach.} After the introduction of the lower
bounds of the Aronson-Caffarelli type, it became clear that the
size of the initial $\LL^1$ or $\LL^p$ norm in a certain ball can be
used in a natural way to define intrinsic quantities for later
times, and this is the approach the authors followed in \cite{BV}
for the easier range $m_c<m<1$\,. The Harnack inequalities we
derive below are based on such an idea and apply also for $0<m\le
m_c$. Indeed, if one wants to apply the result of DiBenedetto et. al.
\cite{DBbook, DGV} mentioned above, to a local weak solution defined
on $[0,T]\times\Omega$, where $T$ is possibly the extinction time,
the Harnack inequality of \cite{DKV,DGV} reads:

\noindent\textit{There exists positive constants $\overline{\delta}<\overline{c}$ depending only on $m,d$ such that if
\begin{equation}\label{Intr.hyp}
\overline{c}\,u(t_0,x_0) \le \left[\frac{\min\{t_0,T-t_0\}}{(8R)^2}\right]^{\frac{1}{1-m}}\qquad\mbox{and }\qquad \dist(x_0,\partial\Omega)<\frac{R}{8}\,,
\end{equation}
we then have that
\[
\overline{c}\,u(t_0,x_0)\le \inf_{x\in B_R(x_0)}u(t,x)\,,
\]
for all times $t_0-\overline{\delta}\,u(t_0,x_0)^{1-m}\,R^2<t<t_0
+\overline{\delta}\,u(t_0,x_0)^{1-m}\,R^2$. The constants
$\overline{\delta}$ and $\overline{c}$ tend to zero as
$m\to 1$ or as $m\to m_c$\,.}
The intrinsic hypothesis \eqref{Intr.hyp} is guaranteed in the good range by the fact that solutions with initial data in $\LL^1_{\rm loc}$ are bounded, while in the very fast diffusion range hypothesis
\eqref{Intr.hyp} fails, and should be replaced by :
\[
u(t,x_0)
    \le \frac{c_{m,d} }{\varepsilon^{\frac{2p\vartheta_p}{1-m}}} \left[\frac{\|u(t_0)\|_{\LL^p({B_R})}\,R^d}{\|u(t_0)\|_{\LL^1({B_R})}R^{\frac{d}{p}}}
            \right]^{2p\vartheta_p}\;\left[\frac{t_0}{R^2}\right]^{\frac{1}{1-m}}\,.
\]
This local upper bound can be derived by the smoothing effect of Theorem \ref{thm.upper}, whenever $t_0+\varepsilon t_*(t_0)<t<t_0+t_*(t_0)$, see full details  in the proof of Theorem \ref{Thm.Parab.Harn3}.

\medskip

\noindent \textsc{The Size of Intrinsic Cylinders.}
We will show that the new critical time
\begin{equation}\label{t*s}
t_*(s)=c_{m,d}\,R^{2-d(1-m)}\|u(s)\|^{1-m}_{\LL^1(B_R(x_0))}
\end{equation}
introduced in Part I, gives the size of the intrinsic cylinders:
in the supercritical fast diffusion range this time can be chosen
a priori just in terms of the initial datum, but in the subcritical
range its size changes with time; roughly speaking the diffusion is
so fast that the initial local information is not relevant after some
time, which is represented by $t^*$. We must bear in mind that a
large class of solutions  completely extinguish in finite time.
We proceed next with the new results.

\medskip

\textbf{Inequalities of Forward, Backward and Elliptic Type.} For
small times cf. Theorem \ref{Thm.Parab.Harn2}, or for suitable
intrinsic cylinders, cf. Theorem \ref{Thm.Parab.Harn3}, we obtain
inequalities where the infimum is taken at a later time than the
supremum (forward Harnack inequalities), or at the same time (elliptic
Harnack inequalities), or even at an earlier time (backward Harnack
inequalities).

Throughout this section we take $0<m<1$ and consider a local nonnegative weak solution $u$
of the FDE defined in a cylinder $Q=(0,T)\times \Omega$, taking
initial data $ u(0,x)=u_0(x)$ in $ \LL^p_{\rm loc}\big(\Omega\big)$\,, with $p=1$ if $m_c<m<1$
or $p>p_c$ if $0<m\le m_c$. We make no assumption on the boundary condition (apart from nonnegativity). Also, let   $x_0$ be a point in $\Omega$ and let $6R \le \dist(x_0,\partial \Omega)$. As  before, we let $T_m$ be the so-called minimal life time, corresponding to data $u_0$ and ball $B_R(x_0)$, and we define $t_*(s)$ as in \eqref{t*s} and $t_*=t_*(0)$, which is equal or less than $T_m$.
First we prove Harnack inequalities for initial times

\begin{thm}\label{Thm.Parab.Harn2}
Under the above conditions, for any $t_0\in(0,t_*]$ and $0\le\vartheta\le\min\{t_*-t_0, t_0/2\}$ the following inequality holds
\begin{equation}\label{Parabolic.Harnack.intr}
\inf_{x\in B_R(x_0)}u\big(t_0\pm\theta ,x\big)
        \ge \mathcal{H}\,u\big(t_0,x_0\big)
\end{equation}
where
\begin{equation}\label{Ht}
\mathcal{H}\,
    =\mathcal{C}_6\,R^{\frac{2-d}{m}}\,
        \left[\frac{\|u_0\|_{\LL^1({B_R})}}{T_m^{\frac{1}{(1-m)}}}\right]^{\frac{1}{m}}
        \left[\frac{\|u_0\|_{\LL^p({B_R})}^{2p\vartheta_p}}{t_0^{{\frac{2p\vartheta_p}{1-m}}}}\,
                    +\frac{1}{R^{\frac{2}{1-m}}}\right]^{-1}
\end{equation}
and $\mathcal{C}_6$ depends only on $m,d,p$\,.  $\mathcal{H}$ goes
to zero when $t_0\to 0$.
\end{thm}

This form of Harnack estimate we propose must be called generalized,
since the constant depends on the solution through certain norms of
the data. But we remind the reader that a proper restriction of the
class of initial data allows to control $\cal H$ at any time $0<t<t_*$.

\noindent {\sl Proof.~} The proof  consists of two steps.

\noindent \textit{From center to infimum.} First of all we have to pass from the center to the minimum in the positivity estimates of Theorem 1.1. Fix a point $z\in \overline{B_R(x_0)}$, ad consider the following MDP, centered at $z$:
\begin{equation}\label{z-sub.Dirichlet.Problem}
\left\{\begin{array}{lll}
 \partial_t u =\Delta (u^m) & ~ {\rm in}~ Q_{T,R_0}=(0,T)\times B_{9R/2}(z)\\[2mm]
 u(0,x)=u_0(x)\chi_{B_R(x_0)} & ~{\rm in}~ B_{R}(x_0),~~~{\rm and}~~~ \supp(u_0)\subseteq B_{2R}(z)\\[2mm]
 u(t,x)=0 & ~{\rm for}~  t >0 ~{\rm and}~ x\in\partial B_{9R/2}(z)\,,
\end{array}\right.
\end{equation}
it is then clear that $B_{R}(x_0)\subset B_{2R}(z)\subset
B_{9R/2}(z)\subset B_{6R}(x_0)$. Applying the result of Theorem \ref{Flux.Thm.Posit.2} to the solution $u$ to the above minimal problem with minimal life time $T_m = T_m(u_0)$, we get
\begin{equation}\label{t*.z}
 t_*\,
    :=\frac{k_0}2\,\left(\frac{9}{2}R-2R\right)^2
        \left[\frac{\int_{B_{2R}(z)}u_0 \dx}{\vol\big(B_{9R/2}(z)\setminus B_{2R}(z)\big)} \right]^{1-m}
     =k_{m,d}\,R^{2-d(1-m)}
     \left[\int_{B_{R}(x_0)}u_0 \dx\right]^{1-m}\le T_m\,,
\end{equation}
which does not depend on $z(t)$. We then obtain:
\[\begin{split}
u^m(t,z)
    &\ge c_1'\,(2R)^{2-d} t^{\frac{m}{1-m}} T_m^{-\frac{1}{1-m}}\;
  \int_{B_{2R}(z)}u_0(x)\dx
  =c_1\,R^{2-d} t^{\frac{m}{1-m}} T_m^{-\frac{1}{1-m}}\; \int_{B_{R}(x_0)}u_0(x)\dx.
\end{split}
\]
for any $t\in [0,t_*]$. Since $z\in B_{R}(x_0)$ is arbitrary and does not enter neither in the above lower bounds, neither in the formula \eqref{t*s} for $t_*=t_*(0)$, we can set $z=z(t)$ as the point such that:
\begin{equation}\label{inf.lower}
\inf_{x\in B_{R}}u(t,x)=u\big(t,z(t)\big)\ge
c_1\,\left[\frac{\int_{B_{R}(x_0)}u_0(x)\dx}{T_m^{\frac{1}{(1-m)}}R^{\,d-2}}\right]^{\frac{1}{m}}
    \,t^{\frac{1}{1-m}}
\end{equation}
for any $0\le t\le t_*$, where $t_*=t_*(0)$ is given by \eqref{t*s}.
Once we have obtained the result for the minimal Dirichlet problem we pass to a general weak solution
as it has been done in Subsection \ref{sect.pls}, hence estimate \eqref{inf.lower} holds for any weak solution.

\noindent\textit{Joining upper and lower estimates.}
Let now $t_0\in(0,t_*]$ and choose $\theta>0$  small so that $t_0/2< t_0-\theta$, and $t_0+\theta\leq t_{*}$. By the lower estimate  \eqref{inf.lower} we know that $u(t_0,x_0)$ is
positive for any $t_0\in(0,t_*]$ and
\begin{equation}\label{low.bfe}
\begin{split}
\inf_{x\in B_R(x_0)}u\big(t_0\pm\theta ,x\big)
    &\ge c_1'^{\frac{1}{m}}\,R^{\frac{2-d}{m}}\|u_0\|_{\LL^1({B_R})}^{\frac{1}{m}}
        T_m^{-\frac{1}{m(1-m)}}\;(t_0\pm\theta)^{\frac{1}{1-m}}\\
    &\ge c_1'' 2^{-\frac{1}{1-m}}\|u_0\|_{\LL^1({B_R})}^{\frac{1}{m}}
        T_m^{-\frac{1}{m(1-m)}}\,R^{\frac{2-d}{m}}\,t_0^{\frac{1}{1-m}}
\end{split}
\end{equation}
We now use the local smoothing effect of Theorem \ref{thm.upper} for $t=t_0\in(0,t_*)$\,:
\begin{equation}\label{sup.final.2}
\begin{split}
u(t_0,x_0)
    &\le\mathcal{C}_3 \left[\frac{\|u_0\|_{\LL^p({B_R})}^{2p\vartheta_p}}{t_0^{\frac{2p\vartheta_p}{1-m}}}\,R^{\frac{2}{1-m}}
      +1\right]\;\left[\frac{t_0}{R^2}\right]^{\frac{1}{1-m}}.
\end{split}
\end{equation}

We have thus proved that for $t_0\in(0,t_*]$ and $\theta\le \min\{t_*-t_0, t_0/2\}$ we have
\begin{equation}\label{Harnack.3}
\inf_{x\in B_R(x_0)}u\big(t_0\pm\theta ,x\big)
    \ge c_2\,\|u_0\|_{\LL^1({B_R})}^{\frac{1}{m}} T_m^{-\frac{1}{m(1-m)}}
         \,R^{\frac{2-d}{m}}\,
      \left[\frac{\|u_0\|_{\LL^p({B_R})}^{2p\vartheta_p}}{t_0^{\frac{2p\vartheta_p}{1-m}}}\,
      +\frac{1}{R^{\frac{2}{1-m}}}\right]^{-1}u(t_0,x_0)\,.
\end{equation}
This concludes the proof.\qed

By shifting the interval $[0,t_*]$ to $[t_0, t_0+t_*(t_0)]\subseteq [0,T]$ we can prove a more intrinsic flavored version of backward-forward-elliptic Harnack inequality in the following
\begin{thm}\label{Thm.Parab.Harn3}
Under the above conditions,  there exists constants $h_1\,,h_2$ depending only on $m,d,p$, such that, for any $\varepsilon\in[0,1]$  the following inequality holds
\begin{equation}\label{Parabolic.Harnack.intr.2}
\inf_{x\in B_R(x_0)}u\big(t\pm \vartheta,x\big)
        \ge h_1\,\varepsilon^{\frac{2p\vartheta_p}{1-m}}\left[\frac{\|u(t_0)\|_{\LL^1({B_R})}R^{\frac{d}{p}}}
            {\|u(t_0)\|_{\LL^p({B_R})}R^d}\right]^{2p\vartheta_p+\frac{1}{m}} u\big(t,x_0\big)
\end{equation}
for any
\[
t_0+\varepsilon t_*(t_0)< t\pm \vartheta < t_0+t_*(t_0)\,,\qquad t_*(t_0)
=h_2\,R^{2-d(1-m)}\|u(t_0)\|^{1-m}_{\LL^1(B_R(x_0))}
\]
\end{thm}

\noindent {\sl Proof of Theorem \ref{Thm.Parab.Harn3}.~}
Assume $t_0=0$, the result will follow by translation in time.
We continue the proof of Theorem \ref{Thm.Parab.Harn2}: we further estimate the
 local smoothing effect of Theorem \ref{thm.upper} for $t=t_0\in(0,t_*)$\,:
\begin{equation}\label{sup.final.2}
\begin{split}
u(t_0,x_0)
    &\le\mathcal{C}_3
        \left[\frac{\|u_0\|_{\LL^p({B_R})}^{2p\vartheta_p}}{t_0^{\frac{2p\vartheta_p}{1-m}}}\,R^{\frac{2}{1-m}}
            +1\right]\;\left[\frac{t_0}{R^2}\right]^{\frac{1}{1-m}}
      \le \mathcal{C}_4 \left[\frac{\|u_0\|_{\LL^p({B_R})}^{2p\vartheta_p}}{\varepsilon
            t_*^{\frac{2p\vartheta_p}{1-m}}}\,R^{\frac{2}{1-m}}
                \right]\;\left[\frac{t_0}{R^2}\right]^{\frac{1}{1-m}}\\
    &\le \frac{\mathcal{C}_5}{\varepsilon^{\frac{2p\vartheta_p}{1-m}}} \left[\frac{\|u_0\|_{\LL^p({B_R})}\,R^d}{\|u_0\|_{\LL^1({B_R})}R^{\frac{d}{p}}}
            \right]^{2p\vartheta_p}\;\left[\frac{t_0}{R^2}\right]^{\frac{1}{1-m}}\\
\end{split}
\end{equation}
since we have put $t_0\ge \varepsilon t_*$ and $t_*=t_*(0)$ as in \eqref{t*s}.
Next we use the estimate for the extinction time proved in
Sections \ref{sec.upper.fet.1} and \ref{sec.upper.fet.2}\,, that can be rewritten as
\[
T_m^{\frac{1}{1-m}}\le k_{m,p,d} R^{\frac{2}{1-m}-\frac{d}{p}}
    \|u_0\|_{\LL^p(B_R)}\qquad\mbox{for any $p\ge \max\{p_c,1\}$}.
\]
The lower estimates \eqref{low.bfe} becomes
\begin{equation}\label{low.bfe2}
\begin{split}
\inf_{x\in B_R(x_0)}u\big(t_0\pm\theta ,x\big)
    &\ge c_1\,\left[\frac{\|u_0\|_{\LL^1({B_R})}R^{\frac{2}{1-m}-d}}
        {T_m^{\frac{1}{(1-m)}}}\right]^{\frac{1}{m}}
        \;\left[\frac{t_0}{R^2}\right]^{\frac{1}{1-m}}
    \ge c_2\,\left[\frac{\|u_0\|_{\LL^1({B_R})}R^{\frac{d}{p}}}
        {\|u_0\|_{\LL^p(B_R)}R^d}\right]^{\frac{1}{m}}
        \;\left[\frac{t_0}{R^2}\right]^{\frac{1}{1-m}}\\
\end{split}
\end{equation}
Joining now inequalities \eqref{sup.final.2} and \eqref{low.bfe2} we obtain that
there exists constants $h_1\,,h_2$ depending only on $m,d,p$, such that, for any $\varepsilon\in[0,1]$  the following inequality holds
\begin{equation}\label{Parabolic.Harnack.intr.2.casi}
\inf_{x\in B_R(x_0)}u\big(t\pm \vartheta,x\big)
        \ge h_1\,\varepsilon^{\frac{2p\vartheta_p}{1-m}}\left[\frac{\|u(0)\|_{\LL^1({B_R})}R^{\frac{d}{p}}}
            {\|u(0)\|_{\LL^p({B_R})}R^d}\right]^{2p\vartheta_p+\frac{1}{m}} u\big(t,x_0\big)
\end{equation}
for any
\[
\varepsilon t_*(0)< t\pm \vartheta < t_*(0)\,,\qquad t_*(0)
=h_2\,R^{2-d(1-m)}\|u(0)\|^{1-m}_{\LL^1(B_R(x_0))}\,.
\]
We conclude the proof by translating the result from $[0, t_*(0)]$ to $[t_0,t_0+t_*(t_0)]$, that we know to be included in $[0,T]$ as explained in Part I\,.\qed

\medskip
\noindent\textbf{Remarks. }
\noindent (i) Estimate \eqref{Parabolic.Harnack.intr.2} is
completely of local type, since it involves only local quantities.
In the supercritical range $m_c<m<1$, we can let $p=1$ in
\eqref{Parabolic.Harnack.intr.2} to get
\begin{equation*}
\inf_{x\in B_R(x_0)}u\big(t\pm \vartheta,x\big)
        \ge h_1\,\varepsilon^{\frac{2p\vartheta_p}{1-m}}u\big(t,x_0\big)
\end{equation*}
for any $t_0+\varepsilon t_*(t_0)< t\pm \vartheta < t_0+t_*(t_0)$.
In this way we recover the above mentioned results of DiBenedetto et al., cf.
\cite{DKV, DGV}. This theorem complements and supports the lower Harnack inequality
\eqref{GoodFDE.Low.Harnack} of Part I, in which the constant does not depend on $u_0$.

Joining the upper  and lower estimates for the Cauchy Problem, we obtain the
Global Harnack principle as the authors did in \cite{BV}.

\noindent(ii) In the subcritical range $0<m\le m_c$, the Harnack estimates cannot have a universal
constant independent of $u_0$, as already mentioned, cf. also
in \cite{DGV} for a counterexample. We have proved that if one allows
the constant to depend on the initial data, then it is possible to
obtain intrinsic Harnack inequalities, and the price  we pay is having
the minimal life time in the constant, as in Theorem \ref{Thm.Parab.Harn2}, but this information is a bit
unpractical and we replace it with some  local $\LL^p$-norm with $p>p_c$ of the initial datum.

\noindent(iii) We have shown that the size of the intrinsic cylinders is always proportional to a
ratio of local $\LL^p$ norms. Note that in the supercritical range it simplifies and only depends on the local $\LL^1$ norm.

\noindent (iv) The quantity $\varepsilon$ represents an arbitrary small waiting time, that is needed in order for the regularization to take place and to allow quantitative intrinsic Harnack inequalities.

\noindent(v) Backward Harnack inequalities are a bit surprising, but they reflect a typical feature of the fast diffusion processes, that is the extinction phenomena, namely
\begin{equation}\label{Parabolic.Harnack.intr.backw}
\inf_{x\in B_R(x_0)}u\big(t- \vartheta,x\big)
        \ge h_1\,\varepsilon^{\frac{2p\vartheta_p}{1-m}}\left[\frac{\|u(t_0)\|_{\LL^1({B_R})}R^{\frac{d}{p}}}
            {\|u(t_0)\|_{\LL^p({B_R})}R^d}\right]^{2p\vartheta_p+\frac{1}{m}} u\big(t,x_0\big)
\end{equation}
for any
\[
t_0+\varepsilon t_*(t_0)< t-\vartheta < t_0+t_*(t_0)\,,\qquad t_*(t_0)
=h_2\,R^{2-d(1-m)}\|u(t_0)\|^{1-m}_{\LL^1(B_R(x_0))}\,.
\]
This inequality is compatible with the fact that the solution extinguish at some later time, remaining
strictly positive before. This backward inequality is typical of singular equation and can not hold for
the degenerate -porous media- case $m>1$, neither for the linear heat equation case, $m=1$.

The same remark applies for the Elliptic Harnack inequality, that is when $\vartheta=0$.

\

\textbf{An alternative form of Harnack Inequalities.}
We provide a form of Harnack inequalities of forward, backward and elliptic type, avoiding the intrinsic framework, and the waiting time $\varepsilon\in[0,1]$.

\begin{thm}\label{Thm.Parab.Harn.0} Under the above conditions, there exists positive constants $\mathcal{C}_{1}$, $\mathcal{C}_{2}$ and $h_2$ depending only on $m,d$ and $p$ such that
\begin{equation}\label{Parabolic.Harnack.alt}
\begin{split}
\sup_{x\in  B_{ R}}u(t,x)
    &\le\frac{\mathcal{C}_{1}}{t^{d\vartheta_p}}\,\|u(t_0)\|_p^{2p\vartheta_p}
     +\mathcal{C}_{2}\left[\frac{\|u(t_0)\|_{\LL^p(B_R)}R^d}
        {\|u(t_0)\|_{\LL^1({B_R})}R^{\frac{d}{p}}}\right]^{\frac{1}{m}}\inf_{x\in B_{R}}u(t\pm\vartheta,x)\\
\end{split}
\end{equation}
for any
\[
0\le t_0<t\pm \vartheta<t_0+t_*(t_0)\le T\,,\qquad t_*(t_0)
=h_2\,R^{2-d(1-m)}\|u(t_0)\|^{1-m}_{\LL^1(B_R(x_0))}\,,
\]
where $\vartheta_p=1/(2p-d(1-m))$.
\end{thm}

\noindent {\sl Proof.} First we observe that we can pass from the center $x_0$ to the supremum in the upper estimate of Theorem \ref{thm.upper} by doubling the radius of the ball on the right hand side, namely

\noindent \textit{There exist positive constants $\mathcal{C}_{1}$, $\mathcal{C}'_{2}$ depending only on $m,d$, such that for any  $t, R>0$  we have}
\begin{equation*}
\begin{split}
\sup_{x\in B_{R}(x_0)}u(t,x)
    &\le\frac{\mathcal{C}_{1}}{t^{d\vartheta_p}}\|u(t_0)\|_{\LL^p(B_{2R}(x_0))}^{2p\vartheta_p}
     +\mathcal{C}'_{2}\left[\frac{t}{R^2}\right]^{\frac{1}{1-m}}.\\
\end{split}
\end{equation*}
Joining the above inequality with the lower bound of Theorem \ref{Thm.Posit} in the form of \eqref{low.bfe2}, we obtain the inequality \eqref{Parabolic.Harnack.alt} for $t_0=0$.
We conclude the proof by shifting the interval $[0, t_*(0)]$ to $[t_0,t_0+t_*(t_0)]$\,,
that we know to be included in $[0,T]$ as explained in Part I\,.\qed

\medskip

\noindent\textbf{Remark.} In the good fast diffusion range we can let $p=1$, so that inequality \eqref{Parabolic.Harnack.alt} reads
\begin{equation*}
\begin{split}
\sup_{x\in  B_{ R}}u(t,x)
    &\le\frac{\mathcal{C}_{1}}{t^{d\vartheta_p}}\,\|u(t_0)\|_1^{2\vartheta_1}
     +\mathcal{C}_{2}\inf_{x\in B_{R}}u(t\pm\vartheta,x)\,.\\
\end{split}
\end{equation*}

\section{\large Concluding section}

First, we sketch a panorama of the local estimates in the different
exponent ranges $m<1$ and for different integrability exponents
$p\ge 1$ that may serve as further orientation for the reader. Then
we make a series of general comments, and finally we review
related works.

\subsection{\large Panorama and open problems}

The values of $m_c$ and $p_c$ are defined in the Introduction.

\begin{figure}[ht]
\centering
\includegraphics[width=13cm,height=5.5cm]{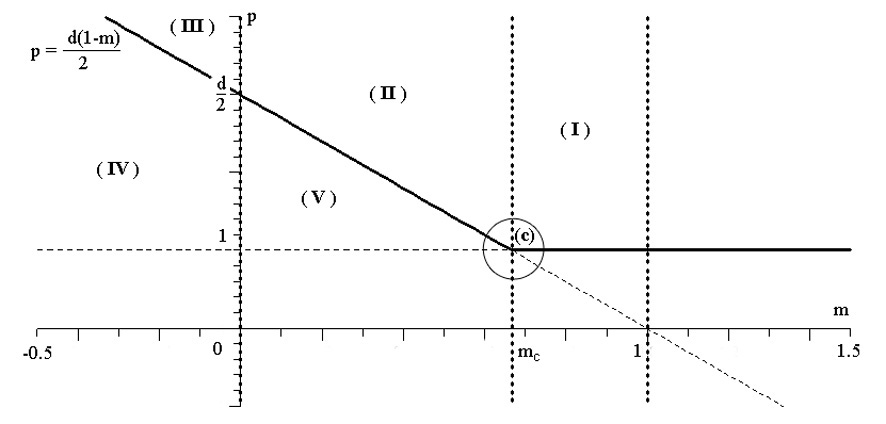}
\end{figure}

\begin{itemize}
\item[(I)] \textsl{Good Fast Diffusion Range:}
$m\in(m_c,1)$ and $p\ge 1$. Here, the local smoothing effect holds,
cf. Theorem  \ref{thm.upper}, as well as  the Reverse smoothing
effect and the global Harnack principle, proved in \cite{BV},
see also \cite{BV3}. As for the result of the present paper,
we have provided a different proof of the positivity result in
Theorem \ref{GoodFDE.Posit}. We have also proved intrinsic Harnack
inequalities of forward, elliptic, or backward type, cf. Theorems
\ref{Thm.Parab.Harn2}, \ref{Thm.Parab.Harn3} and \ref{Thm.Parab.Harn.0}.
recovering the existing results. Finally, for times close to the
extinction time, in case extinction occurs, the authors show in
\cite{BV2, BV3}, that Elliptic Harnack inequalities hold up to extinction.

\item[(II)] \textsl{Very Fast Diffusion Range:} $m\in(0,m_c)$ and $p\ge p_c>1$\,.
The local smoothing effect of Theorem \ref{thm.upper} holds, as well as the
lower estimates of Theorem \ref{Thm.Posit}. These are the only known local
positivity and smoothing results in this range. For any positive time we have
the Aronson-Caffarelli type estimates of Theorem \ref{Thm.Posit}. We have also
proved intrinsic Harnack inequalities of forward, elliptic, or backward type,
cf. Theorems \ref{Thm.Parab.Harn2}, \ref{Thm.Parab.Harn3} and \ref{Thm.Parab.Harn.0}.
These are the only parabolic Harnack inequalities known in this whole range.

\textsl{An open problem is to pass from local to global estimates in this
very fast diffusion range.}  For $m>m_c$ this is done in
\cite{BV} in the form of a global Harnack principle.

\textsl{Another open problem is to find the rate of convergence to an
appropriate extinction profile for the Dirichlet problem on domains.}
A minimal rate of extinction can be obtained by our intrinsic Harnack
estimates, of Theorem \ref{Thm.Parab.Harn3}, details will appear separately.

\item[(c)] \textsl{Critical case:} $m=m_c$ and $p>p_c=1$. The local upper
and lower estimates of zone (II) apply, as well as the Harnack
inequalities. In an upcoming paper we will show how to pass from
local to global estimates, obtaining global lower bounds with
super-exponential time decay.

\item[(III)] \textsl{Negative exponent range:} $m\le 0$ with $p>p_c$.
No positivity result is known in this range, and the technique
developed in this paper does not allow to treat this case. Recall
that solutions of the homogeneous Dirichlet  problem on bounded
domains vanish instantaneously (hence, there is no actual solution).
The local upper bound of Theorem \ref{thm.upper} is still valid, and
this is the only known local upper estimate in this range. \textsl{It is an
open problem in this case to find  positivity and a posteriori Harnack
inequalities, if any.}

\item[(IV)] \textsl{Negative range:} $m\le 0$ with $p<p_c$.
The smoothing effect is not true, since initial data are not in
$\LL^p$ with $p>p_c$, cf. \cite{VazLN}, and the solutions of the
Cauchy problem with data in $\LL^p(\RR^d)$ will vanish
instantaneously, setting a strong limitation to positivity results,
which must be based of course on local bounds. Again, solutions to
the homogeneous Dirichlet  problem on bounded domains vanish
instantaneously. \textsl{It is an open problem to find positivity
estimates in this range, if any.} In general, Harnack inequalities
are not possible in this range since solution may not be (neither
locally) bounded.

\item[(V)] \textsl{Very Fast Diffusion Range,}
$m\in(0,m_c)$ with small integrability exponent $p\in[1,p_c]$. It is
well known that the smoothing effect is not true in general, since
initial data are not in $\LL^p$ with $p>p_c$, cf. \cite{BV, ChVArch,
VazLN}\,. The lower estimates in this case are given by Theorem \ref{Thm.Posit}.
These are the only known local positivity results in this whole range.
In general, Harnack inequalities are not possible in this case since
solution may not be (neither locally) bounded.
\end{itemize}

\subsection{\large Some general remarks}
\label{ss.sgr}

\begin{itemize}
\item We stress  the fact that our results are
completely local, and they apply to any kind of initial-boundary
value problem, in any Euclidean domain: Dirichlet, Neumann, Cauchy,
or problem for large solutions, namely  when $u=+\infty$ on the
boundary, etc. Natural extensions are fast diffusion problems with
variable coefficients and fast diffusion problems on
manifolds.

\item We calculate (almost) explicitly
all the constants, through all the paper.

 \item Our positivity and Harnack inequalities generalize
the results of \cite{DK, DKV, DGV}, valid only in the good fast diffusion
range, in the sense that we recover their result with a different proof,
and we extend it quantitatively to the very fast diffusion range.

\item More specific information can be obtained when we restrict our
attention to particular classes of solutions, like the solutions of
the Cauchy Problem, the solutions of the initial and boundary value
problem on a bounded domain with $u=0$ on the boundary, or with
$u=+\infty$ on the boundary. All these solutions have additional
properties that can be exploited. We will not enter into those
topics for reasons of space.

\item We will not enter either into the derivation of H\"older
continuity and further regularity from the Harnack inequalities.
This is a subject extensively treated in the works of DiBenedetto et
al., see \cite{DiBenedettourbVesp, DBbook, DGV} and references
therein.

\item Actually, the fast diffusion equation can be well-posed even
for Borel measures  (i.e., not locally finite measures) as initial
data. This is done for the Cauchy problem in the whole space in
\cite{ChVArch},  but in that case the smoothness of the constructed
solutions is lost and the concept of extended continuous solution
has to be introduced. The problem turns out to be well-posed in that
extended class in the good range $m_c<m<1$.

 \item   The main ideas of this paper can be extended to related
equations, such as the $p$-Laplacian, but the technical details may not be immediate.
\end{itemize}

\vspace{-.5cm}

\subsection{\large Elliptic connection}\label{Ell.Conn}
Part of the
difficulties of the FDE in the lower range of $m$ can be explained
by the intimate relation of the equation with  semilinear elliptic
theory. Indeed, if we try the separation of variables ansatz
$u(t,x)=H(t)F(x)$ on equation \eqref{FDE} we easily get the
possible formulas $H_1(t)=t^{1/(1-m)}$ or
$H_2(t)=(T-t)^{1/(1-m)}$, and then $G=F^m$ satisfies the elliptic
equation
\begin{equation}\label{ellp.eq.G}
 -\Delta G \pm c\, G^q=0, \quad \mbox{with } \ q=\frac1m\,,
\end{equation}
so that $q>1$ if $0<m<1$. The constant is $c=m/(1-m)$, and the sign
$\pm$ corresponds to the choices $H_1$ or $H_2$ respectively. It is
well-known that the theory of equation (\ref{ellp.eq.G}) is
difficult for large values of $q$, notably for $q\ge
q_s=(d+2)/(d-2)$. The exponent corresponding to $m_c$ is
$q_c=d/(d-2)$, a lower exponent that appears sometimes in the study
of singularities. Note that when $m$ is negative $G$ is an inverse
power of $F$; moreover, $q$ is negative.\footnote{We leave to the
reader the exponential formulas that are obtained for $m=0$.} The
FDE-elliptic connection extends to the study of self-similar
solutions of different types, which is a fundamental tool of the FDE
theory and is described in \cite{Ki93} and \cite{VazLN} among other
references.

It is interesting to check the correspondence of our time evolving
results with the theory of elliptic equations. An easy way of doing
that it to apply the results to special solutions. The simplest
case, i.e., stationary solutions, is too simple, hence we prefer to
try the separate variable solutions
$$
u(t,x)=(T-t)^{1/(1-m)}F(x).
$$
 Putting $U=F^m$ and $q=1/m$ we get the elliptic equation $\Delta U
 + cU^q=0$, as in \eqref{ellp.eq.G}. It can also be written as
$$
\Delta U + a(x)U=0, \quad a(x)=cU^{q-1}=cF^{1-m}.
$$
Let us check the upper estimate \eqref{upper.est} of Theorem
\ref{thm.upper}. Using the separate variable form of $u$ we
immediately see  that the time dependencies disappear (a
confirmation of the correct scaling of the formula) and the we
obtain a local boundedness result for $U$ of the form
$$
U(x)\le C \|U\|_q^{\theta(q)} + C_2R^{-2p/(1-p)},
$$
on the condition that $F\in \LL^p_{loc}(\Omega)$ with $p\ge p_c$ which
means that $a(x)\in \LL^r_{loc}(\Omega)$ with $r\ge d/2$, a classical
condition. This is for us another way of checking that $p_c$ is a
natural exponent.

A similar conclusion can be derived from application of the lower
estimate \eqref{AC.Estimates} of Theorem \ref{Thm.Posit}. We leave
the details to the reader. The elliptic conclusions can also be checked on selfsimilar
solutions of the type $u(t,x)=(T-t)^{\alpha}F(x\,(T-t)^{\beta})$, as
the ones considered in \cite{VazLN}.

\subsection{\large Short review of related works}

The range $m\le m_c$ has remained outside of most of the
publications on the questions of positivity and Harnack estimates.
For positivity and boundedness,  let us first mention the works of
Bertsch and collaborators \cite{BDU90, BU90} who treat the
equation satisfied by the so-called pressure variable
$v=c/u^{1-m}$, i.e., an inverse power of $u$. It covers the
equivalent to the whole fast diffusion range $m<1$ in terms of
viscosity solutions;  the questions are somewhat different from
our program.

We list next and comment on a number of previous results on the
subject of Harnack inequalities  for the fast diffusion equation.
\begin{itemize}

\item DiBenedetto and Kwong proved in \cite{DK} that in the good
fast diffusion range $m_c<m<1$ intrinsic Harnack Inequalities of forward type
do hold, under the positivity assumption that there exists a point
$(t_0,x_0)$ such that $u(t_0,x_0)>0$\,. This value controls from
below the infimum in a small ball at a later time, with sizes
depending on $u(x_0,t_0)$.

\item Later on DiBenedetto, Kwong and Vespri \cite{DKV} improved
on the previous result proving the Global Harnack Principle in the
wider range $m_s=\frac{d-2}{d+2}<m<1$, by means of comparison with
the separation of variable solution, always under the assumption of
positive solutions and a stronger assumption on the initial data,
namely $u_0^m\in W^{1,2}_0(\Omega)$. These estimates are global
in space but not in time since the constants blow up as $t\to 0$.
Hence the interest in combining them with  information we provide
for all small times in direct dependence of the local $\LL^p$ norms
of the initial data.

\item More recently, DiBenedetto, Gianazza and Vespri, \cite{DGV},
extended the results to the variable coefficient case in the form
Harnack inequalities which are of Forward, Elliptic and Backward
type; it applies in the good FDE range, and always under the
positivity assumption  for some $(x_0,t_0)$. Some of these estimates
had been proved by the authors in the constant coefficient case in
\cite{BV}\,, see also \cite{BV3}.

\item  The above mentioned Harnack inequalities imply H\"older
continuity of the solution and sometimes analyticity, cf. \cite{DKV, DGV}.

\item The power $m_s=(d-2)/(d+2)$, has been studied by
Del Pino and Saez \cite{DPS}, as part of the study of the
asymptotics of the evolutionary Yamabe  problem. They perform the
transformation into a fast diffusion problem posed on the sphere via
stereographic projection, which is possible for this exponent. They
get an elliptic Harnack inequality which holds for a good class of
solutions, but they do not prove a parabolic Harnack inequality.

\item None of the above quoted papers considers the problem of finding
Harnack inequalities when the time approaches the finite extinction
time (if there is one). This has been done by the authors in
\cite{BV2, BV3}, showing that Elliptic Harnack inequalities hold up
to the extinction time. The proof is completely different, we draw
fine asymptotic properties, by a careful analysis of the extinction
profile.
\end{itemize}

Summing up, two results are known in the lower range $m\le m_c$:
\cite{DPS} that applies for $m=m_s$, and \cite{DKV} that applies for
$m_s<m<1$. They both refer to a different point of view.

\medskip

\section*{\large Appendixes}
\setcounter{section}{5}

\subsection*{\large A1\; Aleksandrov's Reflection Principle}

Here we state the Reflection Principle of Aleksandrov in a slightly
different form, more useful to our purposes. We already used this proposition,
in \cite{BV}. Other forms of the same principle, in
different settings can be found, for example in \cite{GV}),
Proposition 2.24 (pg. 51) or in \cite{Ar-Pe}, Lemma 2.2.

\begin{prop}[Local Aleksandrov's Reflection Principle,
\cite{BV}\,]\label{Local.Aleks}
\noindent Let $B_{\lambda R_0}(x_0)\subset\RR^d$ be an open ball
with center in $x_0\in\RR^d$ of radius $\lambda\,R_0$ with $R_0>0$
and $\lambda>2$. Let $u$ be a solution to problem
\begin{equation}\label{FDE.Problem.Aleks}
\begin{split}
\left\{\begin{array}{lll}
\partial_t u =\Delta (u^m) & ~ {\rm in}~ (0,+\infty)\times B_{\lambda R_0}(x_0)\\
u(0,x)=u_0(x) & ~{\rm in}~ B_{\lambda R_0}(x_0) \\
u(t,x)=0 & ~{\rm for}~  t >0 ~{\rm and}~ x\in\partial B_{\lambda R_0}(x_0)\\
\end{array}\right.
\end{split}
\end{equation}
with $\supp(u_0)\subset B_{R_0}(x_0)$. Then, for any $t>0$ one has:
\[
u(t,x_0)\ge u(t,x_2)
\]
for any $t>0$ and for any $x_2\in A_{\lambda, R_0}(x_0)=B_{\lambda
R_0}(x_0)\setminus B_{2R_0}(x_0)$. Hence,
\begin{equation}\label{Aleks.Mean}
u(t,x_0)\ge \left|A_{\lambda,
R_0}(x_0)\right|^{-1}\int_{A_{\lambda,
R_0}(x_0)}u(t,x)\dx=\oint_{A_{\lambda, R_0}(x_0)}u(t,x)\dx
\end{equation}
\end{prop}
\noindent The proof of this result can be found in the Appendix of \cite{BV}\,.

\subsection*{\large A2\; On the extinction time}

The extinction time plays a role in our study of positivity  in Part
I. We devote some paragraphs to comment on its occurrence in FDE. It
has been observed in the literature, cf. \cite{BGV-JEE, BV, VazLN},
that Lemma \ref{Lem.Local.L1.norms} can be used to obtain lower
bounds on the finite extinction time $T=T(u_0)$. Indeed, just let
$t=0$, $s=T$ and $v\equiv 0$, in \eqref{Lem.Local.Lp.p=1}, to get
\begin{equation}\label{Low.Bound.Ext.Time}
c_1(m,d,\lambda)R^{d(1-m)-2}
\|u_0\|_{\LL^1(B_R)}^{1-m}\le \;T
\end{equation}
Suppose that $\Omega=\RR^d$. When $m_c<m<1$, it is easy to see
that $2-d(1-m)>0$, so that, simply by letting $R\to\infty$, we see
that $T(u_0)\to+\infty$ if $u_0\in\LL^1(\RR^d)$, that means that
solutions corresponding to initial data $u_0\in\LL^1(\RR^d)$ do
not extinguish in finite time, i.e., there is global positivity.

Such a situation does not occur for solutions defined in all of
$\RR^d$ when $0<m<m_c$. Indeed, in that fast diffusion range
solutions with initial data in $\LL^1(\RR^d)$ may extinguish. The
question is studied in Chapters 5--7 of \cite{VazLN} where upper
and lower estimates for $T$ in terms of the data are obtained for
the problem posed in the whole space $\RR^d$ with $m<m_c$.

On the other hand, solutions extinguish in finite time for the
Cauchy-Dirichlet problem posed in a bounded domain with zero
boundary conditions for all $0<m<1$, cf. subsections \ref{sec.upper.fet.1},
\ref{sec.upper.fet.2} or \cite{BGV, DiDi}.

An  estimate similar to \eqref{Low.Bound.Ext.Time} cannot  be valid
for $m< 0$ since it is known that the Cauchy problem does not admit
solutions with data $u_0\in \LL^1(\RR^d)$ \cite{VazNonex}, or even
in $\LL^p(\RR^d)$ with $p<p_c$ \cite{DasDP} for $m\le 0$. We may say
in these cases that $T=0$. Actually, when one tries to solve
approximate problems with, say, strictly positive and bounded
boundary data and pass to the limit, the approximate solutions
converge to zero uniformly in cylinders of the form
$Q_\tau=(\tau,\infty)\times \RR^d$. This can be summed up as the
formation of an initial discontinuity layer, cf. \cite{VazLN}. No
solutions exist either for the Cauchy-Dirichlet problem posed in a
bounded domain with zero boundary conditions when $m\le 0$. There is
also a peculiar case,  $m=0$ and $d=2$, where there exist solutions
but the waiting time can be fixed a priori independently of the
initial data \cite{VazLN}, hence no  estimate as above is possible
either.

Note finally that the $\LL^p$ estimate of Theorem \ref{Lem.Local.Lp.norms}
for  $p>1$ cannot be used to obtain lower estimates on FET, since they
hold only for $t\ge s\ge 0$.

\subsection*{\large A3\; Details of the iterative calculations of Subsection \ref{sect3}, Step 3}

We show here in detail of a calculation used to pass to the limit
when $k\to \infty$ in the inequality \eqref{Iter.step.k+1}\,, in the
proof of Theorem \ref{Thm.1.Integrated}. We adopt the notations used
there.
\[
\begin{split}
\prod_{j=1}^{k} j^{2\frac{(1+\frac{1}{q})^{k-j}}{p_{k+1}}}
    &\le \exp\left[ 2\sum_{j=1}^{k}\frac{(1+\frac{1}{q})^{k-j}}{p_{k+1}}\log(j)\right]
    \le \exp\left[ 2\sum_{j=1}^{k}\frac{(1+\frac{1}{q})^{k-j}\log(j)}
        {\big[p_0-q(1-m)\big]\left[1+\frac{1}{q}\right]^{k+1}+(q+1)(1-m)}\right]\\
    &\le s_0 \exp\left[ 2\sum_{j=1}^{k}(1+\frac{1}{q})^{-j-1}\log(j)\right]
    \le s_1 \exp\left[ 2\sum_{j=1}^{k}(1+\frac{1}{q})^{-j}\right]\\
        &\le s_1 \exp\left[ 2\sum_{j=1}^{\infty}(1+\frac{1}{q})^{-j}\right]
        =s_1\ee^{2(q+1)}\\
\end{split}
\]
where $s_i$, $i=0,1$ are positive numerical constant.\qed

\

\vspace{1cm}

\noindent {\textbf{\large \sc Acknowledgment.} Both authors funded
by Project MTM2005-08760-C02-01 (Spain) and European Science
Foundation Programme ``Global and Geometric Aspects of Nonlinear
Partial Differential Equations''. The first author supported by
Juan de la Cierva (Spain) Programme.

\vskip 1cm



\end{document}